\documentclass[a4paper]{article}
\usepackage{amssymb}
\usepackage{amsmath}
\usepackage{amsthm}
\usepackage{authblk}
\usepackage[USenglish]{babel}
\usepackage{esint}
\usepackage{graphicx}
\usepackage[colorlinks=true]{hyperref}
\usepackage[latin1]{inputenc}
\newtheorem{theorem}{Theorem}[section]
\newtheorem{lemma}[theorem]{Lemma}
\newtheorem{proposition}[theorem]{Proposition}
\newtheorem{corollary}[theorem]{Corollary}
\newtheorem{remark}[theorem]{Remark}
\newtheorem{definition}[theorem]{Definition}

\newcommand{\N}{\mathbb N}
\newcommand{\Z}{\mathbb Z}

\newcommand{\R}{\mathbb R}
\newcommand{\C}{\mathbb C}
\renewcommand{\S}{\mathbb S}
\renewcommand{\P}{\mathbb P}
\renewcommand{\H}{\mathbb H}
\newcommand{\B}{\mathbb B}

\newcommand{\mbb}{\mathbb}
\newcommand{\mbf}{\mathbf}
\newcommand{\mcal}{\mathcal}
\newcommand{\mfrak}{\mathfrak}

\newcommand{\mrm}{\mathrm}

\renewcommand{\a}{\alpha}
\renewcommand{\b}{\beta}

\newcommand{\G}{\Gamma}
\renewcommand{\d}{\delta}
\newcommand{\D}{\Delta}
\newcommand{\e}{\varepsilon}
\newcommand{\z}{\zeta}
\renewcommand{\t}{\theta}

\newcommand{\la}{\lambda}

\newcommand{\s}{\sigma}
\newcommand{\si}{\varsigma}
\newcommand{\Si}{\Sigma}

\newcommand{\ph}{\varphi}
\renewcommand{\o}{\omega}
\renewcommand{\O}{\Omega}
\newcommand{\wt}{\widetilde}
\newcommand{\wh}{\widehat}
\newcommand{\ol}{\overline}
\newcommand{\ul}{\underline}

\newcommand{\ub}{\underbrace}
\newcommand{\fr}{\frac}
\newcommand{\pa}{\partial}
\newcommand{\n}{\nabla}
\newcommand{\fa}{\forall}
\newcommand{\ex}{\exists}
\newcommand{\es}{\emptyset}
\newcommand{\wk}{\rightharpoonup}

\newcommand{\us}{\underset}

\newcommand{\lot}{\longleftarrow}
\newcommand{\lto}{\longrightarrow}

\newcommand{\sm}{\setminus}
\renewcommand{\Cup}{\bigcup}

\newcommand{\sub}{\subset}
\newcommand{\Sub}{\Subset}

\newcommand{\nin}{\not\in}
\newcommand{\eq}{\equiv}

\newcommand{\pl}{\oplus}

\newcommand{\x}{\times}
\renewcommand{\c}{\circ}
\newcommand{\cd}{\cdot}
\newcommand{\ds}{\dots}

\newcommand{\tx}{\text}
\newcommand{\q}{\quad}
\renewcommand{\l}{\left}
\renewcommand{\r}{\right}
\newcommand{\tbf}{\textbf}
\newcommand{\bthm}{\begin{theorem}}
\newcommand{\ethm}{\end{theorem}}
\newcommand{\blem}{\begin{lemma}}
\newcommand{\elem}{\end{lemma}}
\newcommand{\bprop}{\begin{proposition}}
\newcommand{\eprop}{\end{proposition}}
\newcommand{\bcor}{\begin{corollary}}
\newcommand{\ecor}{\end{corollary}}
\newcommand{\bdefi}{\begin{definition}}
\newcommand{\edefi}{\end{definition}}
\newcommand{\bpf}{\begin{proof}}
\newcommand{\epf}{\end{proof}}
\newcommand{\bl}{\begin{array}{l}}
\newcommand{\bll}{\begin{array}{ll}}
\newcommand{\barr}{\begin{array}}
\newcommand{\earr}{\end{array}}
\newcommand{\bite}{\begin{itemize}}
\newcommand{\eite}{\end{itemize}}
\newcommand{\bequ}{\begin{equation}}
\newcommand{\eequ}{\end{equation}}
\newcommand{\beqa}{\begin{eqnarray}}
\newcommand{\eeqa}{\end{eqnarray}}
\newcommand{\beqy}{\begin{eqnarray*}}
\newcommand{\eeqy}{\end{eqnarray*}}

\allowdisplaybreaks

\begin{document}

\everymath{\displaystyle}

\title{Existence and non-existence results for the $SU(3)$ singular Toda system on compact surfaces}
\author{Luca Battaglia\thanks{Universit\'e Catholique de Louvain, Institut de Recherche en Math\'ematique et Physique, Chemin du Cyclotron $2$, $1348$ Louvain-la-Neuve - luca.battaglia@uclouvain.be},
Andrea Malchiodi\thanks{Scuola Normale Superiore, Piazza dei Cavalieri $7$, $56126$ Pisa (Italy) - andrea.malchiodi@sns.it\\Authors are supported by the PRIN project \emph{Variational and perturbative aspects of nonlinear differential problems}, grant number $201274$FYK$7$.}}
\date{}

\maketitle

\begin{abstract}
\noindent We consider the $SU(3)$ singular Toda system on a compact surface $(\Si,g)$
$$\l\{\bl-\D u_1=2\rho_1\l(\fr{h_1e^{u_1}}{\int_\Si h_1e^{u_1}\mrm dV_g}-1\r)-\rho_2\l(\fr{h_2e^{u_2}}{\int_\Si h_2e^{u_2}\mrm dV_g}-1\r)-4\pi\sum_{m=1}^M\a_{1m}\l(\d_{p_m}-1\r)\\-\D u_2=2\rho_2\l(\fr{h_2e^{u_2}}{\int_\Si h_2e^{u_2}\mrm dV_g}-1\r)-\rho_1\l(\fr{h_1e^{u_1}}{\int_\Si h_1e^{u_1}\mrm dV_g}-1\r)-4\pi\sum_{m=1}^M\a_{2m}\l(\d_{p_m}-1\r)\earr\r.,$$
where $h_i$ are smooth positive functions on $\Si$, $\rho_i\in\R_+$, $p_m\in\Si$ and $\a_{im}>-1$.\\
We give both existence and non-existence results under some conditions on the parameters $\rho_i$ and $\a_{im}$. Existence results are obtained using variational methods, which involve a geometric inequality of new type; non-existence results are obtained using blow-up analysis and localized Poho\v zaev-type identities.
\end{abstract}

\section{Introduction}

Let $\Si$ be a closed surface and $g$ be a Riemannian metric on $\Si$. Consider the following system on $\Si$:
\bequ\label{todaN}
-\D u_i=\sum_{j=1}^2a_{ij}\rho_j(h_je^{u_j}-1)-4\pi\sum_{m=1}^M\a_{im}(\d_{p_m}-1),\q\q\q i=1,2,
\eequ
where $\D=\D_g$ is the Laplace-Beltrami operator, $\rho_i$ are positive parameters, $h_i$ are smooth positive functions on $\Si$, $\a_{im}$ are real numbers greater than $-1$, $p_m$ are given points of $\Si$, and $A=(a_{ij})_{ij}$ is the Cartan matrix of $SU(3)$
$$\l(\barr{cc}2&-1\\-1&2\earr\r).$$
System \eqref{todaN} is known as the $SU(3)$ singular Toda system. Together with its $N\x N$ extension, it has been widely studied in literature due to its important role in both geometry and mathematical physics. In geometry, it appears in the description of holomorphic curves in $\C\P^3$ (see e.g. \cite{cal,cw,bw}), while in mathematical physics it arises in the non-abelian Chern-Simons theory (see \cite{dunne,yang,tar08}). The singularities represent respectively the ramification points of the complex curves and the \emph{vortices} of the wave functions.\\

To better understand this system, it is convenient to re-write it in an equivalent form. Let $G_p$ be the Green function of $-\D$ centered at a point $p\in\Si$, namely the solution of
$$\l\{\bl-\D G_p=\d_p-1\\\int_\Si G_p\mrm dV_g=0\earr\r..$$
Consider now the substitution
$u_i\mapsto u_i+4\pi\sum_{m=1}^M\a_{im}G_{p_m}$: the newly-defined $u=(u_1,u_2)$ solves
$$\l\{\bl-\D u_1=2\rho_1\l(\wt h_1e^{u_1}-1\r)-\rho_2\l(\wt h_2e^{u_2}-1\r)\\-\D u_2=2\rho_2\l(\wt h_2e^{u_2}-1\r)-\rho_1\l(\wt h_1e^{u_1}-1\r)\earr\r.,$$
where the new functions $\wt h_i$ have the expression
$$\wt h_i:=h_ie^{-4\pi\sum_{m=1}^M\a_{im}G_{p_m}}\q\q\q i=1,2$$
and verify
\bequ\label{hi}
\wt h_i\in C^\infty\l(\Si\sm\{p_1,\ds,p_M\}\r)\q\q\q \wt h_i|_{\Si\sm\{p_1,\ds,p_M\}}>0\q\q\q\wt h_i\sim d(\cd,p_m)^{2\a_{im}}\q\tx{near }p_m.
\eequ
Integrating by parts over the whole $\Si$, we deduce
$$\int_\Si\wt h_1e^{u_1}\mrm dV_g=\int_\Si\wt h_1e^{u_2}\mrm dV_g=1,$$
therefore the system is equivalent to
\bequ\label{toda}
\l\{\bl-\D u_1=2\rho_1\l(\fr{\wt h_1e^{u_1}}{\int_\Si\wt h_1e^{u_1}\mrm dV_g}-1\r)-\rho_2\l(\fr{\wt h_2e^{u_2}}{\int_\Si\wt h_2e^{u_2}\mrm dV_g}-1\r)\\-\D u_2=2\rho_2\l(\fr{\wt h_2e^{u_2}}{\int_\Si\wt h_2e^{u_2}\mrm dV_g}-1\r)-\rho_1\l(\fr{\wt h_1e^{u_1}}{\int_\Si\wt h_1e^{u_1}\mrm dV_g}-1\r)\earr\r.,
\eequ

Problem \eqref{toda} admits a variational formulation, that is its solutions are critical points of the following energy functional defined on $H^1(\Si)^2$:
\bequ\label{jrho}
J_\rho(u):=\int_\Si Q(u)\mrm dV_g-\sum_{i=1}^2\rho_i\l(\log\int_\Si\wt h_ie^{u_i}\mrm dV_g-\int_\Si u_i\mrm dV_g\r).
\eequ
Here, $Q(u)$ is given by
$$Q(u)=\fr{|\n u_1|^2+\n u_1\cd\n u_2+|\n u_2|^2}3,$$
$\n=\n_g$ is the gradient given by the metric $g$ and $\cd$ denotes the Riemannian scalar product.\\
To study the properties of the functional $J_\rho$, a basic tool is the Moser-Trudinger inequality, which was proved in \cite{batmal,bat2} (and, for the regular case, in \cite{jw01}).
\bequ\label{mttoda}
4\pi\sum_{i=1}^2\min\l\{1,1+\min_m\a_{im}\r\}\l(\log\int_\Si\wt h_ie^{u_i}\mrm dV_g-\int_\Si u_i\mrm dV_g\r)\le\int_\Si Q(u)\mrm dV_g+C.
\eequ
As a consequence, $J_\rho$ is bounded from below as long as $\rho_i\le4\pi\min\l\{1,1+\min_m\a_{im}\r\}$ for both $i=1,2$. Moreover, if both parameters are strictly smaller than these thresholds, the functional is coercive in the space of functions with zero average; there will be no loss of generality in restricting the problem to this space, since both \eqref{toda} and \eqref{jrho} are invariant by addition of constants. Hence, in this case we get minimizing solutions.\\

If one (or both) of the $\rho_i$ is allowed to attain greater values, then one can build suitable test functions to show that the energy functional is unbounded from below, as was done in the same papers where \eqref{mttoda} is proved. Therefore, one can no longer use minimization techniques to find critical points. However it is possible to prove that when the Euler-Lagrange energy \eqref{jrho} becomes largely negative at least one of the functions $\wt h_i e^{u_i}$ has to concentrate near a finite number of points. One can eventually derive existence results out of this statement using min-max or Morse theory.\\

To describe in more detail the situation we first consider \emph{Liouville's equation}, that is the scalar counterpart of \eqref{toda}:
$$-\D u=2\rho\l(\fr{he^u}{\int_\Si he^u\mrm dV_g}-1\r)-4\pi\sum_{m=1}^M\a_m(\d_{p_m}-1).$$
Through a change of variable similar to that before \eqref{toda}, this is equivalent to
\bequ\label{liou}
-\D u=2\rho\l(\fr{\wt he^u}{\int_\Si\wt he^u\mrm dV_g}-1\r).
\eequ
with $\wt h$ having the same behavior as in \eqref{hi} around singular points.\\
Liouville's equation has also great importance in geometry and mathematical physics: it appears in the problem of prescribing the Gaussian curvature on surfaces with conical singularities and in models from abelian Chern-Simons theory. This problem has also been very much studied in literature, with many results concerning existence of solutions, compactness properties, blow-up analysis et al., which have been summarized e.g. in the reviews \cite{mal10,tar10}.\\

\eqref{liou} is the Euler-Lagrange equation for the functional 
\bequ\label{irho}
I_\rho(u):=\fr{1}2\int_\Si|\n u|^2\mrm dV_g-2\rho\l(\log\int_\Si\wt he^u\mrm dV_g-\int_\Si u\mrm dV_g\r).
\eequ
The classical Moser-Trudinger inequality and its extension to the singular cases 
(\cite{mos,fon,chenwx,tro}) yield boundedness from below of $I_\rho$ if and only if $\rho\le4\pi\min\l\{1,1+\min_m\a_m\r\}$ and coercivity if and only if $\rho$ is strictly smaller than this value.\\

For larger values of $\rho$, despite the lack of lower bounds on the energy $I_\rho$, it is however possible to prove that functions with low energy must \emph{concentrate} near finitely-many points. A heuristic reason for this fact goes as follows: the Moser-Trudinger inequality can be localized on any region of $\Sigma$ via cut-off functions, see \cite{cl91}. A consequence of this fact is that functions that are \emph{spread} over $\Si$ satisfy a Moser-Trudinger inequality with an improved constant, which favors lower bounds on $I_\rho$. Hence, if lower bounds fail, $u$ should \emph{concentrate} rather than spread.\\
Notice that when all the $\a_j$'s are negative the localized Moser-Trudinger constant near a singular point $p_i$ is $4\pi(1+\a_i)$, while near a regular point it is simply $4\pi$. Based on these considerations, in \cite{carmal} the following weighted cardinality $\o_{\ul\a}$ on finite sets was introduced:
\bequ\label{wbar}
\o_{\ul\a}(\{x\}):=\l\{\bll1+\a_m&\tx{if }x=p_m\\1&\tx{if }x\nin\{p_1,\ds,p_M\}\earr\r.\q\q\q\o_{\ul\a}\l(\Cup_k\{x_k\}\r):=\sum_k\o_{\ul\a}(\{x_k\});
\eequ
and it was shown that if a function $u$ has low energy, then the normalized measure $\wt he^u$ must distributionally approach the following set of measures (appeared also in \cite{dm})
$$\Si_{\rho,\ul\a}:=\l\{\sum_{x_k\in\mcal J}t_k\d_{x_k}:\,x_k\in\Si,\,t_k\ge0,\,\sum_{x_k\in\mcal J}t_k=1,\,4\pi\o_{\ul\a}(\mcal J)<\rho\r\}.$$
Using variational methods, a compactness result in \cite{bt1} and a monotonicity argument in \cite{str} it was also shown that, endowing $\Si_{\rho,\ul\a}$ with the weak topology of distributions, solutions to \eqref{liou} (up to a discrete set of $\rho$'s, for compactness reasons) exist provided $\Si_{\rho,\ul\a}$ is non-contractible. We notice that the problem is not always solvable, as in the classical case of the \emph{teardrop}: the sphere with only one singular point. Sufficient and necessary conditions for contractibility were given in \cite{car}.\\
The case of positive singularities was treated in \cite{bdm} on surfaces with positive genus. There are some other existence results (\cite{barmal,mr11}) which also work for the case of the sphere or of the real projective plane. We also refer to \cite{cl13} for the derivation of a degree-counting formula.\\

We turn now to system \eqref{toda}: for the \emph{regular case} some existence results were found in \cite{mn} ($\rho_1<4\pi$ and $\rho_2\nin4\pi\N$), in \cite{mr13} ($\rho_i\in(4\pi,8\pi)$), \cite{jkm} ($\rho_1\in(4\pi,8\pi)$ and $\rho_2\nin4\pi\N$) and in \cite{bjmr} ($\Si$ of positive genus and $\rho_i\nin4\pi\N$). In the latter paper, with a construction related to that in \cite{bdm}, the case of positive singular weights was also treated while in \cite{bat1}, still for positive genus, some cases with negative coefficients $\a_{im}$ were discussed.\\

The above reasoning for the scalar singular equation allows to prove a related alternative for the two components of the system. If we use the compact notation $\ul\a_1:=(\a_{11},\ds,\a_{1m}),\,\ul\a_2:=(\a_{21},\ds,\a_{2m})$, then it turns out that for $J_\rho(u)$ low either $\wt h_1e^{u_1}$ is distributionally close to $\Si_{\rho_1,\ul\a_1}$ or $\wt h_2e^{u_2}$ is close to $\Si_{\rho_2,\ul\a_2}$. To express this (non-exclusive) alternative, it is natural introduce the \emph{join} of two topological spaces $X$ and $Y$ (see for instance \cite{hat}):
\bequ\label{join}
X\star Y:=\fr{X\x Y\x[0,1]}\sim,
\eequ
where $\sim$ is the equivalence relation among triples $(x,y,t)$ given by
$$(x,y,0)\sim(x,y',0)\q\fa\,x\in X,\,\fa\,y,y'\in Y\q\q\q(x,y,1)\sim(x',y,1)\q\fa\,x,x'\in X,\,\fa\,y\in Y.$$
The join of $\Si_{\rho_1,\ul\a_1}$ and $\Si_{\rho_2,\ul\a_2}$ could then be used to characterize low-energy levels of $J_\rho$, with the join parameter $s\in[0,1]$ expressing whether $\wt h_1e^{u_1}$ is closer to $\Si_{\rho_1,\ul\a_1}$ or $\wt h_2e^{u_2}$ is closer to $\Si_{\rho_2,\ul\a_2}$ (for example $s=\fr{1}2$ would describe couples with the same scale of concentration).\\
This description is however not optimal in general, as it does not take accurately into account the interaction between two components $u_1$ and $u_2$. For the regular case of \eqref{toda}, in \cite{mr13} it was shown that the \emph{relative rate of concentration} of the two components plays a role in this matter.\\
More precisely, it was shown that if $u_1,u_2$ concentrate near the same point and with the same scale (see Section \ref{s:prel} for a more precise definition of the latter), then the Moser-Trudinger constants for the system double. As a consequence of this fact it turns out that, when $\rho_1,\rho_2\in(4\pi,8\pi)$ and no singularities occur, then join elements of the form $\l(x,x,\fr{1}2\r)$, $x\in\Si$ have to be excluded (see \cite{jkm} for higher values of $\rho_1$).\\

One of the main goals of this paper is to show a new improved inequality for the singular system \eqref{toda}, in order to understand at the same time the effect of the interaction of the two components among themselves and with the singularities. We prove in particular (see Section \ref{s:impr}) that if the two components are concentrated near the same singular point with the same rate, coercivity of the Euler-Lagrange energy holds provided $\rho_1,\rho_2<4\pi(2+\a_1+\a_2)$ (notice that with no extra assumption coercivity holds under the weaker condition $\rho_i<4\pi(1+\a_i)$ for all $i$'s).\\

We expect these new improved inequalities would allow us to prove existence results in rather general cases. However for simplicity here we restrict ourselves to relatively low values of $\rho_1,\rho_2$, in such a way that the above-defined measures $\Si_{\rho_i,\ul\a_i}$ are supported in at most one singular point of $\Si$. Precisely, defining the two numbers
\bequ\label{rhobar}
\ol\rho_1:=4\pi\min\l\{1,\min_{m\ne m'}(2+\a_{1m}+\a_{1m'})\r\}\q\q\q\ol\rho_2:=4\pi\min\l\{1,\min_{m\ne m'}(2+\a_{2m}+\a_{2m'})\r\},
\eequ
by choosing $\rho_i<\ol\rho_i$, $\Si_{\rho_i,\ul\a_i}$ will contain only Dirac deltas centered at singular points $p_m$ for some $m\in\{1,\ds,M\}$.
In fact, $\rho_i<4\pi$ regular points are excluded, while $\rho_i<\ol\rho_i$ ensures the one-point support condition.\\

The first main result contained in this paper is the following one. We would need to exclude some null set $\Gamma$ of $\R^2$ for compactness reasons, see Section \ref{s:prel}.\\

\bthm\label{ex}
Let $\G$ be as in \eqref{gamma}, $(\ol\rho_1,\ol\rho_2)$ be as in \eqref{rhobar}, and let $\rho\in\R^2_+\sm\G$ satisfy $\rho_i<\ol\rho_i$ for both $i=1,2$.\\
Define integer numbers $M_1,M_2,M_3$ by:
\beqa
\label{m123}
\nonumber M_1:=\#\{m:\,4\pi(1+\a_{1m})<\rho_1\}\q\q\q M_2:=\#\{m:\,4\pi(1+\a_{2m})<\rho_2\}\\
M_3:=\#\{m:\,4\pi(1+\a_{im})<\rho_i\tx{ and }\rho_i<4\pi(2+\a_{1m}+\a_{2m})\tx{ for both }i=1,2\}.
\eeqa
Then system \eqref{toda} admits solutions provided the following condition holds 
$$(M_1,M_2,M_3)\nin\{(1,m,0),(m,1,0),(2,2,1),(2,3,2),(3,2,2),\,m\in\N\}.$$
\ethm

\begin{remark}
We will see  that the above assumptions on the $M_i$'s are necessary: in fact we will get a non-existence 
result for every case not covered by the theorem. 
\end{remark}

By the previous description low sub-levels of $J_\rho$ can be identified with the topological join of $\Si_{\rho_1,\ul\a_1}$ and of $\Si_{\rho_2,\ul\a_2}$, with some points removed. Under the assumptions on the $\rho_i$'s this join consists of a graph $\mcal{X}$ made of segments whose end-points belong to $\{p_1,\dots,p_m\}$. For a more precise description of it we refer to Section \ref{s:top}, where some pictures are also included.  The conditions on $(M_1, M_2, M_3)$ 
in the previous theorem ensure that this graph is non-contractible.\\

The second part of this paper, see Section \ref{s:nonex} will be devoted to the proof of some non-existence results, showing that in general 
some assumptions on the parameters $\rho_i$ are necessary to get existence of solutions. We begin by considering a simple situation: the unit disk of $\R^2$ with a singularity at the origin, and solutions satisfying Dirichlet boundary conditions.\\

\bthm\label{nonex1}
Let $\l(\mbb B^2,g_0\r)$ be the standard unit disk, suppose $h_1,h_2\eq1,\,M=1$ and let $\a_1,\a_2>-1$ be the 
singular weights of the point $p=0\in\mbb B$.
If $\rho$ satisfies
$$\rho_1^2-\rho_1\rho_1+\rho_2^2-4\pi(1+\a_1)\rho_1-4\pi(1+\a_2)\rho_2\ge0,$$
then there are no solutions to the system
\bequ\label{todadisk}
\l\{\bl-\D u_1=2\rho_1\fr{|x|^{2\a_1}e^{u_1}}{\int_\B|x|^{2\a_1}e^{u_1(x)}\mrm dx}-\rho_2\fr{|x|^{2\a_2}e^{u_2}}{\int_\B|x|^{2\a_2}e^{u_2(x)}\mrm dx}\\-\D u_2=2\rho_2\fr{|x|^{2\a_2}e^{u_2}}{\int_\B|x|^{2\a_2}e^{u_2(x)}\mrm dx}-\rho_1\fr{|x|^{2\a_1}e^{u_1}}{\int_\B|x|^{2\a_1}e^{u_1(x)}\mrm dx}\\u_1|_{\pa\B}=u_2|_{\pa\B}=0\earr\r..
\eequ
\ethm

This result is proved via the Poho\v zaev identity, and extends a scalar one from \cite{barmal}.\\
With a similar proof, one can find non-existence for \eqref{toda} on the standard sphere with one singular point or two antipodal ones. We remark that, as for Theorem \ref{nonex1}, the following result still holds if we allow the coefficients $\a_{im}$ to be positive, thus showing that the general existence result contained in \cite{bjmr} cannot be extended to spheres.\\
As shown by pictures in Section \ref{s:nonex}, non-existence occurs on a region delimited by four curves: we get two or three connected components, which intersect the axis $\rho_1$ in the segment joining $(0,4\pi(1+\a_{21}))$ and $(0,4\pi(1+\a_{22}))$, thus including the scalar case considered in \cite{barmal}. Moreover, such regions also include some  cases which are not covered by Theorem \ref{ex}, in particular $(M_1,M_2,M_3)\in\{(1,m,0),(m,1,0),(2,2,1)\}$.\\

\bthm\label{nonex2}
Let $(\Si,g)=\l(\S^2,g_0\r)$ be the standard sphere, suppose $h_1,h_2\eq1,\,M=2$, let $(\a_{11},\a_{21})\ne(\a_{12},\a_{22})$ be the weights of the antipodal points $\{p_1,p_2\}\sub\S^2$, with 
$\a_{im}>-1$. If either
\bequ\label{cond1}
\l\{\bl\rho_1^2+\rho_2^2-\rho_1\rho_2-4\pi(1+\a_{11})\rho_1-4\pi(1+\a_{21})\rho_2\le0\\
\rho_1^2+\rho_2^2-\rho_1\rho_2-4\pi(1+\a_{12})\rho_1-4\pi(1+\a_{22})\rho_2\ge0\\
\rho_1^2-\rho_2^2-4\pi(1+\a_{11})\rho_1+4\pi(1+\a_{22})\rho_2\le0\\
\rho_1^2-\rho_2^2-4\pi(1+\a_{12})\rho_1+4\pi(1+\a_{21})\rho_2\ge0\earr\r. 
\eequ
and at least one inequality is strict, or if all the opposite inequalities hold, then system \eqref{toda} admits no solutions.
\ethm

The third result we present makes no assumptions on the topology of $\Si$. In fact, its proof will use a \emph{localized} blow-up analysis around one singular point, similarly to some result in \cite{car}. We argue by contradiction, assuming that a solution $u^n$ of \eqref{toda} exists for a sequence $(\a_{11}^n,\a_{12}^n)\us{n\to+\infty}\lto(-1,-1)$.\\
Such a sequence must blow-up, hence we consider all the possibilities given by concentration-compactness theorems (from \cite{ln,jw01,batmal}, which we will recall in Section \ref{s:prel}). We will exclude all of these cases but the blow-up around the point $p_1$. Finally, we will also rule this out by a local version of the Poho\v zaev identity, hence getting a contradiction.\\

Just like Theorem \ref{nonex2}, the following result shows the sharpness of assuming \emph{all} the singularities to be non-negative in \cite{bjmr}. In fact, the statement still holds true if we allow all the coefficients $\a_{12},\ds,\a_{1M},\a_{22},\ds,\a_{2M}$ to be positive and only $\a_{11},\a_{12}<0$.\\

\bthm\label{nonex3}
Let $\G_{\ul\a_{1\wh1},\ul\a_{2\wh1}}\sub\mbb R_+^2$ be as in \eqref{gamma}, with $\ul\a_{i\wh1}:=(\a_{i2},\ds,\a_{iM})$ and let $\rho\in\R_+^2\sm\G_{\ul\a_{1\wh1},\ul\a_{2\wh1}}$ and $\a_{12},\ds,\a_{1M},\a_{22},\ds,\a_{2M}$ be fixed.
Then, there exists $\a_*\in(-1,0)$ such that the system \eqref{toda} is not solvable for $\a_{11},\a_{12}\le\a_*$. Moreover, $\a_*$ can be chosen uniformly for $\rho$ in a given $\mcal K\Sub\R_+^2\sm\G_{\ul\a_{1\wh1},\ul\a_{2\wh1}}$.
\ethm

The last non-existence result gives a counterexample to Theorem \ref{ex}, in the case $(M_1,M_2,M_3)=(2,3,2)$ (which was not covered by Theorem \ref{nonex2}). We basically combine arguments from Theorems \ref{nonex2} and \ref{nonex3}: we consider the standard unit sphere, take $\rho_i,\a_{im}$ so that we have $(M_1,M_2,M_3)=(2,3,2)$ and we let one of the parameters $\a_{im}$ go to $-1$. By a blow-up analysis we reduce ourselves to 
 the scalar version of Theorem \ref{nonex2} (see \cite{barmal}, Proposition $5.8$) and we prove that no solution can exist if that coefficient is too close to $-1$.\\

\bthm\label{nonex4}
Let $(\Si,g)=\l(\S^2,g_0\r)$ be the standard sphere, $\ol\rho_1,\ol\rho_2$ be as in \eqref{rhobar}, suppose $h_1,h_2\eq1,\,M=3$ with $p_1,p_2$ being anti-podal and $\a_{13}=0$, and let $\rho_1,\a_{11},\a_{12},\a_{21},\a_{22}$ be fixed so that
$$4\pi(1+\a_{1m})<\rho_1<\max\l\{\ol\rho_2,4\pi(2+\a_{1m}+\a_{2m})\r\}\q\q\q\fa\,m=1,2.$$
Then, there exists $\a_*\in(-1,0)$ such that \eqref{toda} is not solvable if $\a_{23}\le\a_*$ and $\rho_2$ satisfies
$$4\pi(1+\a_{2m})<\rho_2<\max\l\{\ol\rho_2,4\pi(2+\a_{1m}+\a_{2m})\r\}\q\q\q\fa\,m=1,2,3.$$
\ethm

The paper is organized as follows. In Section \ref{s:prel} we provide some notation and preliminary results that will be used later on. In Section \ref{s:top} we introduce the above-mentioned space $\mcal X$ and study its topology and homology groups. Section \ref{s:test} is devoted to the construction of test functions from $\mcal X$ to arbitrarily low sub-levels of $J_\rho$, whereas in Section \ref{s:impr} we prove new improved Moser-Trudinger inequalities which will be used. In Section \ref{s:ex} Theorem \ref{ex} will be proved using the strategy described before. Finally, Section \ref{s:nonex} will be concerned with the non-existence Theorems \ref{nonex1}, \ref{nonex2}, \ref{nonex3} and \ref{nonex4}.\\

\section{Notation and preliminaries}\label{s:prel}

In this section we will provide some notation and some known preliminary results that will be used throughout the rest of the paper.\\

\subsection{Notation}\label{ss:not}

We will denote the indicator function of a set $\O\sub\Si$ as
$$\mbf1_\O(x):=\l\{\bll1&\tx{if }x\in\O\\0&\tx{if }x\nin\O\earr\r..$$
The metric distance between two points $p,q\in\Si$ will be denoted by $d(x,y)$; similarly, for any $\O,\O'\sub\Si$ we will write:
$$d(x,\O):=\inf\{d(x,y):\,x\in\O\}\q\q\q d(\O,\O'):=\inf\{d(x,y):\,x\in\O,\,y\in\O'\}.$$
If $\O$ has a smooth boundary, given $x\in\pa\O$, the outer normal at $x$ will be denoted as $\nu(x)$. The open metric disk centered at $p\in\Si$ with radius $r>0$ will be indicated as $B_r(p)$. For $r_2>r_1>0$ we denote the open annulus centered at $p$ with radii $r_1,r_2$ as
$$A_{r_1,r_2}(p):=\{x\in\Si:\,r_1<d(x,p)<r_2\}=B_{r_2}(p)\sm\ol{B_{r_1}(p)}.$$
If $\O\sub\Si$ has a smooth boundary, for any $x\in\pa\O$ we will denote the outer normal at $x$ as $\nu(x)$.\\
For a given $u\in L^1(\O)$ and a measurable set $\O\sub\Si$ with positive measure, the average of $u$ on $\O$ will be denoted as
$$\fint_\O u\mrm dV_g=\fr{1}{|\O|}\int_\O u\mrm dV_g.$$
In particular, since we are assuming $|\Si|=1$
$$\int_\Si u\mrm dV_g=\fint_\Si u\mrm dV_g.$$
The subset of the space $H^1(\Si)$ consisting of functions with null average is denoted as
$$\ol H^1(\Si):=\l\{u\in H^1(\Si):\,\int_\Si u\mrm dV_g=0\r\}.$$
As recalled before, both the system \eqref{toda} and its energy functional $J_\rho$ defined in \eqref{jrho} are invariant by adding constants to the components $u_i$. Therefore, there will be no loss of generality in restricting our study of the problem on $\ol H^1(\Si)^2$.\\
The sub-levels of $J_\rho$, which, as anticipated, will play an essential role throughout the whole paper, will be denoted as
$$J_\rho^a=\l\{u\in H^1(\Si)^2:\,J_\rho(u)\le a\r\}.$$

We will denote with the symbol $X\simeq Y$ a homotopy equivalence between two topological spaces $X$ and $Y$.\\
The composition of two homotopy equivalences $F_1:X\x[0,1]\to Y$ and $F_2:Y\x[0,1]\to Z$ satisfying $F_1(\cd,1)=F_2(\cd,0)$ is the map $F_2\ast F_1:X\x[0,1]\to Z$ defined by
$$F_2\ast F_1:(x,s)\mapsto\l\{\bll F_1(x,2s)&\tx{if }s\le\fr{1}2\\F_2(x,2s-1)&\tx{if }s>\fr{1}2\earr\r..$$
The identity map on $X$ will be denoted as $\mrm{Id}_X$.\\
$H_q(X)$ will stand for the $q^\mrm{th}$ homology group with coefficient in $\Z$ of a topological space $X$ as $H_q(X)$. An isomorphism between two homology groups will be denoted just by equality sign. Reduced homology groups will be denoted as $\wt H_q(X)$, namely
$$H_0(X)=\wt H_0(X)\pl\Z\q\q\q H_q(X)=\wt H_q(X)\q\tx{if }q\ge1.$$
The $q^\mrm{th}$ Betti number of $X$, namely the dimension of its $q^\mrm{th}$ group of homology, will be indicated by $b_q(X):=\mrm{rank}(H_q(X))$. The symbol $\wt b_q(X)$ will stand for the dimension of $\wt H_q(X)$, that is
$$\wt b_0(X)=b_0(X)-1\q\q\q\wt b_q(X)=b_q(X)\q\tx{if }q\ge1.$$
Throughout the paper we will use the letter $C$ to denote large constants which can vary between different formulas or lines. To stress the dependence on some parameter(s) we may add subscripts such as $C_\a$. We will denote by the symbol $o_\a(1)$ a quantity tending to $0$ as $\a\to0$ or as $\a\to+\infty$. Subscripts will be omitted when they are evident by the context. Similarly, we will use the symbol $x(\a)\sim_\a y(\a)$ to express that the ratio between $x(\a)$ and $y(\a)$ is bounded both from above and below by two positive constants as $\a$ goes to $0$ or to $+\infty$. In other words, $\log\fr{x(\a)}{y(\a)}=O_\a(1)$.\\

\subsection{Compactness results}\label{ss:cpt}

We first state the compactness result for solutions of \eqref{toda}: it can be deduced by a concentration-compactness alternative from \cite{batmal,batman,ln} and a quantization of local blow-up limits from \cite{jlw,lwzhang,wz}.\\
A global compactness result was already given in \cite{batman} using the quantization result from \cite{lwzhang}. In the same way, we here deduce an improvement using \cite{wz}. We present the concentration-compactness theorem in a slightly more general form, which will be useful in the proof of Theorem \ref{nonex3}.\\

\bthm(\cite{ln}, Theorem $4.2$; \cite{batmal}, Theorem $2.1$; \cite{batman}, Theorem $2.1$)\label{conccomp}
Let $\O\sub\Si$ be an open domain and $\{u^n=(u_1^n,u_2^n)\}_{n\in\N}$ be a sequence of solutions of \eqref{toda} on $\O$ with $h_i^n\us{n\to+\infty}\lto h_i>0$ in $C^1\l(\ol\O\r)$ and $\rho_i^n\us{n\to+\infty}\lto\rho_i$. Define
$$\mcal S_i:=\l\{x\in\Si:\,\ex\,\l\{x^n\r\}_{n\in\N}\sub\Si:\,u_i^n(x^n)-\log\int_\Si\wt h_i^ne^{u_i^n}\mrm dV_g\us{n\to+\infty}\lto+\infty\r\}.$$
Then, up to subsequences, one of the following alternatives occurs:
\begin{itemize}
\item(Compactness) For each $i=1,2$ either $u_i^n-\log\int_\Si\wt h_i^ne^{u_i^n}\mrm dV_g$ is uniformly bounded in $L^\infty_{\mrm{loc}}(\O)$ or it tends locally uniformly to $-\infty$.
\item(Blow-up) The blow-up set $\mcal S:=\mcal S_1\cup\mcal S_2$ is non-empty and finite.
\end{itemize}
Moreover,
$$\rho_i^n\fr{\wt h_i^ne^{u_i^n}}{\int_\Si\wt h_i^ne^{u_i^n}\mrm dV_g}\us{n\to+\infty}\wk r_i+\sum_{x\in\mcal S_i}\s_i(x)\d_x$$
in the sense of measures, with $r_i\in L^1(\O)$ and $\s_i(x)$ defined by
$$\s_i(x):=\lim_{r\to0}\lim_{n\to+\infty}\rho_i^n\fr{\int_{B_r(x)}\wt h_i^ne^{u_i^n}\mrm dV_g}{\int_\Si\wt h_i^ne^{u_i^n}\mrm dV_g}.$$
Finally, if $x\in\mcal S_i\sm\{p_1,\ds,p_M\}$ and $2\s_i(x)-\s_{3-i}(x)\ge4\pi$, then $r_i\eq0$. The same holds if $p_m\in\mcal S_i$ and $2\s_i(p_m)-\s_{3-i}(p_m)\ge4\pi(1+\a_{im})$.
\ethm

We next have the following quantization result for $(\s_1(x),\s_2(x))$.

\bthm(\cite{jlw}, Proposition $2.4$; \cite{lwzhang}, Theorem $1.1$; \cite{wz})\label{quant}
Let $\mcal S,\s_i(x)$ be defined as in Theorem \ref{conccomp} and suppose $x\in\mcal S$.
If $x\nin\{p_1,\ds,p_M\}$, then $(\s_1(x),\s_2(x))$ is one of the following:
$$(4\pi,0)\q\q\q(0,4\pi)\q\q\q(4\pi,8\pi)\q\q\q(8\pi,4\pi)\q\q\q(8\pi,8\pi).$$
If $x=p_m$ and $\a_{1m},\a_{2m}<0$, then $(\s_1(p_m),\s_2(p_m))$ is one of the following:
$$(4\pi(1+\a_{1m}),0)\q\q\q(0,4\pi(1+\a_{2m}))\q\q\q(4\pi(1+\a_{1m}),4\pi(2+\a_{1m}+\a_{2m}))$$
$$(4\pi(2+\a_{1m}+\a_{2m}),4\pi(1+\a_{2m}))\q\q\q(4\pi(2+\a_{1m}+\a_{2m}),4\pi(2+\a_{1m}+\a_{2m})).$$
In particular, either $r_1\eq0$ or $r_2\eq0$.
\ethm

By putting together Theorems \ref{conccomp} (applied with $h_i^n=h_i$ on the whole $\Si$) and \ref{quant} we get the following:\\

\bcor\label{global}
Let $\G'_{i,\mcal M}\sub\R_+$ be defined, for $i=1,2$ and $\mcal M\sub\{1,\ds,M\}$, by
$$\G'_{i,\mcal M}:=4\pi\l\{n+\sum_{m'\in\mcal M'}(1+\a_{im'})+\sum_{m\in\mcal M}(2+\a_{1m}+\a_{2m}):\,n\in\N,\,\mcal M'\sub\{1,\ds,M\}\sm\mcal M\r\}.$$
and define $\G=\G_{\ul\a_1,\ul\a_2}$, where 
\bequ\label{gamma}
\G_{\ul\a_1,\ul\a_2}=\Cup_{\mcal M\sub\{1,\ds,M\}}\l(\G'_{1,\mcal M}\x\l[\sum_{m\in\mcal M}4\pi(1+\a_{2m}),+\infty\r)\cup\l[\sum_{m\in\mcal M}4\pi(1+\a_{1m}),+\infty\r)\x\G'_{2,\mcal M}\r).
\eequ
Then the family of solutions $\{u_\rho\}_{\rho\in\mcal K}\sub\ol H^1(\Si)^2$ of \eqref{toda} is uniformly bounded in $W^{2,q}(\Si)^2$ for some $q>1$ for any given $\mcal K\Sub\R_+^2\sm\G$
\ecor

Actually, Theorem \ref{quant} holds in this form only assuming $\a_{1m},\a_{2m}\le\a_0$ for some $\a_0>0$. For general values of $\a_{im}$ a finite number of other local blow-up limits is allowed (see \cite{lwzhang}, Proposition $2.4$ for details), therefore a global compactness result similar to \ref{global} still holds true.
Anyway, all the cases which are not considered in the previously stated results verify $\sigma_i(p_m)\ge4\pi$ for both $i$'s, so as long as we are assuming $\rho_1,\rho_2<4\pi$ the values we have to exclude are all contained in $\G$.\\

Concerning compactness, we have a useful result which can be deduced from minor modifications of the argument in \cite{luc}. It basically states the existence of bounded Palais-Smale sequences for $\rho$ belonging to a dense set of $\R_+^2\sm\G$. Putting together with the compactness result stated before, we get:

\blem\label{deform}
Let $\rho\nin\G$ be given and let $a<b$ be such that \eqref{toda} has no solutions in $\{J_\rho\in[a,b]\}$.
Then, $J_\rho^a$ is a deformation retract of $J_\rho^b$.
\elem

We also deduce that $J_\rho$ is uniformly bounded from above on solutions, hence we have:\\

\bcor\label{contr}
Let $\rho\nin\G$ be given. Then, there exists $L>0$ such that $J_\rho^L$ is a deformation retract of $H^1(\Si)^2$; in particular, it is contractible.
\ecor

From now on, we will always assume to take $\rho\in\R^2_+\sm\G$, except in Section \ref{s:nonex}.\\

\subsection{Moser-Trudinger inequalities and their improved versions}\label{ss:mt}

We have the following Moser-Trudinger inequalities for the scalar Liouville equation and for the Toda system respectively.\\

\bthm(\cite{mos}, Theorem $2$; \cite{fon}, Theorem $1.7$; \cite{chenwx}, Theorem $I$; \cite{tro}, Corollary $10$.)
Let $\wt h$ be as in \eqref{hi}.
Then, there exists $C=C_\Si>0$ such that any $u\in H^1(\Si)^2$ satisfies
\bequ\label{mtscal}
16\pi\min\l\{1,1+\min_m\a_m\r\}\l(\log\int_\Si\wt he^u\mrm dV_g-\int_\Si u_i\mrm dV_g\r)\le\int_\Si|\n u|^2\mrm dV_g+C.
\eequ
Equivalently, $I_\rho$ defined by \eqref{irho} is bounded from below if and only if $\rho\le4\pi\min\l\{1,1+\min_m\a_m\r\}$ and it is coercive if and only if $\rho<16\pi\min\l\{1,1+\min_m\a_m\r\}$.\\
In the latter case, it admits a global minimizer $u$ which solves \eqref{liou}.
\ethm

\bthm(\cite{jw01}, Theorem $1.3$; \cite{batmal}, Theorem $1.1$.)\label{t:mttoda}
Inequality \eqref{mttoda} holds for any $u=(u_1,u_2)\in H^1(\Si)^2$.
Equivalently, $J_\rho$ defined by \eqref{jrho} is bounded from below in $H^1(\Si)^2$ if and only if $\rho_i\le4\pi\min\l\{1,1+\min_m\a_{im}\r\}$ for both $i=1,2$, and it is coercive if and only if $\rho_i<4\pi\min\l\{1,1+\min_m\a_{im}\r\}$.
In the latter case, it admits a global minimizer $u=(u_1,u_2)$ which solves \eqref{toda}.
\ethm

We also need a Moser-Trudinger inequality on manifolds with boundary, which extends the \emph{scalar} inequality from \cite{cy88}.\\
Before the statement, we introduce a class of smooth open subset of $\Si$ which satisfy an \emph{exterior and interior sphere condition} with radius $\d>0$:
\bequ\label{eq:ud}
\mfrak A_\d:=\l\{\O\sub\Si:\,\fa\,x\in\pa\O\,\ex\,x'\in\O,\,x''\in\Si\sm\ol\O:\,x=\ol{B_\d(x')}\cap\pa\O=\ol{B_\d(x'')}\cap\pa\O\r\}
\eequ

\bthm\label{mtdisk}
Take $B:=B_1(0)\sub\R^2$ and $u=(u_1,u_2)\in H^1(B)^2$.
Then, there exists $C>0$ such that
$$2\pi\sum_{i=1}^2\l(\log\int_Be^{u_i(x)}\mrm dx-\fint_Bu_i(x)\mrm dx\r)\le\int_BQ(u(x))\mrm dx+C.$$
The same result holds if $B$ is replaced by a simply connected domain belonging to $\mfrak A_\d$ for some $\d>0$, with the constant $C$ is replaced with some $C_\d>0$.
\ethm

As a sketch of a proof, consider a conformal diffeomorphism from $B$ to the unit upper half-sphere and reflect the image of $u$ through the equator. Now, apply the Moser-Trudinger inequality to the reflected $u'$, which is defined on $\S^2$. The Dirichlet integral of $u'$ will be twice the one of $u$ on $B$, while the average and the integral of $e^{u'}$ will be the same, up to the conformal factor. Therefore the constant $4\pi$ is halved to $2\pi$. Starting from a simply connected domain, one can exploit the Riemann mapping theorem to map it conformally on the unit disk and repeat the same argument. The exterior and interior sphere condition ensures the boundedness of the conformal factor.

From the inequality in Theorem \ref{mtdisk} one can easily deduce a \emph{localized} Moser-Trudinger inequality, arguing via cut-off and Fourier decomposition as in \cite{mr13}.\\

\blem\label{mtloc}
For any $\e>0,\,\a_1,\a_2\in(-1,0]$ there exists $C=C_\e$ such that for any $u\in H^1(B)^2$
\beqa
\label{loc}4\pi\sum_{i=1}^2(1+\a_i)\l(\log\int_{B_\fr{1}2(0)}|x|^{2\a_i}e^{u_i(x)}\mrm dx-\fint_Bu_i(x)\mrm dx\r)\le(1+\e)\int_BQ(u(x))\mrm dx+C,\\
\nonumber\\
\nonumber4\pi(1+\a_1)\l(\log\int_{B_\fr{1}8(0)}|x|^{2\a_1}e^{u_1(x)}\mrm dx-\fint_Bu_1(x)\mrm dx\r)+2\pi\l(\log\int_{A_{\fr{1}4,1}(0)}e^{u_2(x)}\mrm dx-\fint_Bu_2(x)\mrm dx\r)&\le\\
\label{locmix}(1+\e)\int_BQ(u(x))\mrm dx+C.
\eeqa
\elem

We will now discuss some inequalities of \emph{improved type}, which hold for special classes of functions.
First, we will provide a \emph{macroscopic} improved Moser-Trudinger inequality for the Toda system. Basically, if $u_1$ and $u_2$ are spread in different sets at a positive distance within each other, then we can get a better constant than in Theorem \ref{t:mttoda}.
Before stating the improved inequality, let us introduce the space of positive normalized $L^1$ functions
\bequ\label{a}
\mcal A:=\l\{f\in L^1(\Si):\,f>0\tx{ a.e. and }\int_\Si f\mrm dV_g=1\r\}.
\eequ
We can associate to any function $u\in H^1(\Si)^2$ a couple of elements of $\mcal A$, through the map
\bequ\label{fui}
(u_1,u_2)\mapsto\l(\fr{\wt h_1e^{u_1}}{\int_\Si\wt h_1e^{u_1}\mrm dV_g},\fr{\wt h_2e^{u_2}}{\int_\Si\wt h_2e^{u_2}\mrm dV_g}\r)=:(f_{1,u},f_{2,u}).
\eequ
Such a map is easily seen to be continuous, through the scalar Moser-Trudinger inequality \eqref{mtscal}.\\

\blem\label{mtmacro}(\cite{bat1}, Lemma $4.3$)\\
Let $\d>0,J_1,K_1,J_2,K_2\in\N$ be given, let $\{m_{11},\ds,m_{1J_1},m_{21},\ds,m_{2J_2}\}\sub\{1,\ds,M\}$ be a selection of indices, $\{\O_{ij}\}_{i=1,2}^{j=1,\ds,J_i+K_i}$ be measurable subsets of $\Si$ such that
$$\a_{im_{ij}}\le0\q\q\q\fa\,i=1,2,j=1,\ds,J_i$$
$$d(\O_{ij},\O_{ij'})\ge\d\q\q\q\fa\,i=1,2,\,\fa\,j,j'=1,\ds,J_i+K_i,\,j\ne j'$$
$$d(p_m,\O_{ij})\ge\d\q\q\q\fa\,i=1,2,\fa\,j=1,\ds,K_i+M_i,\,\fa\,m=1,\ds,M,\,\,m\ne m_{ij};$$
and $u\in H^1(\Si)^2$ satisfy
$$\int_{\O_{ij}}f_{i,u}\mrm dV_g\ge\d\q\q\q\fa\,i=1,2,\,\fa\,j=1,\ds,J_i+K_i.$$
Then, for any $\e>0$ there exists $C=C_{\Si,\d,J_1,K_1,J_2,K_2,\e}>0$ such that
$$4\pi\sum_{i=1}^2\l(K_i+\sum_{j=1}^{J_i}\l(1+\a_{im_{ij}}\r)\r)\l(\log\int_\Si\wt h_ie^{u_i}\mrm dV_g-\int_\Si u_i\mrm dV_g\r)\le(1+\e)\int_\Si Q(u)\mrm dV_g+C.$$
\elem

Let us recall the weighted barycenters defined in \eqref{wbar}.
These are a subset of the space $\mcal M(\Si)$ of the Radon measures of $\Si$, endowed with the $\mrm{Lip}'$ norm, using duality with Lipschitz functions:
\bequ\label{lip}
\|\mu\|_{\mrm{Lip}'(\Si)}:=\sup_{\phi\in\mrm{Lip}(\Si),\|\phi\|_{\mrm{Lip}(\Si)}\le1}\l|\int_\Si\phi\mrm d\mu\r|.
\eequ

We will denote the distance induced by this norm by $d_{\mrm{Lip}'(\Si)}$. One can easily see that $\mcal M(\Si)$ contains the space $\mcal A$ defined in \eqref{a}.\\
From now on we will assume, until Section \ref{s:ex}, that $\rho_1<\ol\rho_1,\rho_2<\ol\rho_2$ (see \eqref{rhobar}), hence each measure in $\Si_{\rho_i,\ul\a_i}$ will be supported at only one point of $\Si$. Therefore we can identify, with a little abuse of notation, $\d_x\in\Si_{\rho_i,\ul\a_i}$ with $x\in\Si$ and write
$$\Si_{\rho_i,\ul\a_i}:=\l\{x\in\Si:\,4\pi\o_{\ul\a_i}(\{x\})<\rho_i\r\}=\{p_m:\,4\pi(1+\a_{im})<\rho_i\}\sub\Si.$$
Notice that, by choosing in \eqref{lip}, $\phi=d(\cd,y)$ we have $d_{\mrm{Lip}'(\Si)}(\d_x,\d_y)\sim d(x,y)$ for any $x,y\in\Si$. This means that, allowing $\rho$ to attain higher values (as was done in \cite{bat1,carmal}), we get a space which contains a homeomorphical copy of $\Si$.\\

In terms of $\Si_{\rho_i,\ul\a_i}$, from Lemma \ref{mtmacro} we deduce that at least one between $f_{1,u}$ and $f_{2,u}$ is arbitrarily close to the respective weighted barycentric space.\\

We need to define, for each $f\in\mcal A$, a \emph{center of mass} and a \emph{scale of concentration}, inspired by \cite{mr13} (Proposition $3.1$) but such that the center of mass belongs to a given finite set $\mcal F\sub\Si$ (which will be, in our applications, a subset of the singular points). As in \cite{mr13}, we will map $\mcal A$ on the topological cone over $\mcal F$ of height $\d$, which is defined by
\bequ\label{cone}
\mcal C_\d\mcal F:=\fr{\mcal F\x[0,\d]}\sim,
\eequ
where the equivalence relation $\sim$ is given by $(x,\d)\sim(x',\d)$ for any $x\in\Si$. The meaning of such an identification is the following: if a function $f\in\mcal A$ does not concentrate around any point $x\in\mcal F$, then we cannot define a center of mass: in this case we set the scale equals to $\d$, that is \emph{large}.\\

\blem\label{center}
Let $\mcal F:=\{x_1,\ds,x_K\}\sub\Si$ be a given finite set and $\mcal A,\,\mcal C_\d$ be defined by \eqref{a} and \eqref{cone}. Then, for $\d>0$ small enough there exists a map $\psi=(\b,\si)=(\b_\mcal F,\si_\mcal F):\mcal A\to\mcal C_\d\mcal F$ such that:
\bite
\item If $\si(f)=\d$, then either $\int_{\Si\sm\Cup_{x\in\mcal F}B_\d(x)}f\mrm dV_g\ge\d$ or there exists $x',x''\in\mcal F$ with $x'\ne x''$ and
$$\int_{B_\d(x')}f\mrm dV_g\ge\d\q\q\q\int_{B_\d(x'')}f\mrm dV_g\ge\d$$
\item If $\si(f)<\d$, then
$$\int_{B_{\si(f)}(\b(f))}f\mrm dV_g\ge\d\q\q\q\int_{\Si\sm B_{\si(f)}(\b(f))}f\mrm dV_g\ge\d.$$
\eite
Moreover, if $f^n\us{n\to+\infty}\lto\d_x$ for some $x\in\mcal F$, then $(\b(f^n),\si(f^n))\us{n\to+\infty}\lto(x,0)$.
\elem

\bpf
Fix $\tau\in\l(\fr{1}2,1\r)$, take $\d\le\fr{\min_{x,x'\in\mcal F,\,x\ne x'}d(x,x')}2$ and define, for $k=1,\ds,K$,
$$I_k(f):=\int_{B_\d(x_k)}f\mrm dV_g;\q\q\q\q\q\q I_0(f):=\int_{\Si\sm\Cup_{x\in\mcal F}B_\d(x)}f\mrm dV_g=1-\sum_{k=1}^kI_k(f),$$
Choose now indices $\wt k, \wh k$ such that
$$I_{\wt k}(f):=\max_{k\in\{0,\ds,K\}}I_k(f)\q\q\q\q\q\q I_{\wh k}(f):=\max_{k\ne\wt k}I_k(f).$$
We will define the map $\psi$ depending on $\wt k$ and $I_{\wt k}(f)$:
\bite
\item $\wt k=0$. Since $f$ has little mass around each of the points $x_k$, we set $\si(f)=\d$ and do not define $\b(f)$, as it would be irrelevant by the equivalence relation in \eqref{cone}. The assertion of the lemma is verified, up to taking a smaller $\d$, because
$$\int_{\Si\sm\Cup_{x\in\mcal F}B_\d(x)}f\mrm dV_g=I_0(f)\ge\fr{1}{K+1}\ge\d$$
\item $\wt k\ge1,\,I_{\wt k}(f)\le\fr{K\tau}{1-\tau}I_{\wh k}(f)$. Here, $f$ has still little mass around the point $x_{\wt k}$ (which could not be uniquely defined), so again we set $\si(f):=\d$. It is easy to see that $I_{\hat{k}}(f)\ge\fr{1-\tau}K$, so
$$\int_{B_\d\l(x_{\wt k}\r)}f\mrm dV_g\ge\int_{B_\d\l(x_{\wh k}\r)}f\mrm dV_g\ge\fr{1-\tau}K$$
\item $\wt k\ge1,\,I_{\wt k}(f)\ge\fr{K\tau}{1-\tau}I_{\wh k}(f)$. Now, $I_{\wt k}(f)>\tau$, so one can define a scale of concentration $s\l(x_{\wt k},f\r)\in(0,\d)$ of $f$ around $x_{\wt k}$, uniquely determined by
$$\int_{B_{s\l(x_{\wt k},f\r)}\l(x_{\wt k}\r)}f\mrm dV_g=\tau.$$
We can also define a center of mass $\b(f)=x_{\wt k}$ but we have to interpolate for the scale:
\bite
\item Case $I_{\wt k}(f)\le\fr{2K\tau}{1-\tau}I_{\wh k}(f)$: setting
$$\si(f)=s\l(x_{\wt k},f\r)+\fr{I_{\wt k}(f)}{\fr{K\tau}{1-\tau}I_{\wh k}(f)}\l(\d-s\l(x_{\wt k},f\r)\r),$$
we get $s\l(x_{\wt k},f\r)<\si(f)<\d$; moreover, $I_{\wh k}(f)\ge\fr{1-\tau}{K(1+\tau)}$, hence
$$\int_{B_{\si(f)}(\b(f))}f\mrm dV_g\ge\int_{B_{s\l(x_{\wt k},f\r)}\l(x_{\wt k}\r)}f\mrm dV_g=\tau\ge\d$$
$$\int_{\Si\sm B_{\si(f)}(\b(f))}f\mrm dV_g\ge\int_{\Si\sm B_\d\l(x_{\wt k}\r)}f\mrm dV_g\ge\fr{1-\tau}{K(1+\tau)}\ge\d$$
\item Case $I_{\wt k}(f)\ge\fr{2K\tau}{1-\tau}I_{\wh k}(f)$: we just set $\si(f):s\l(x_{\wt k},f\r)$ and we get
$$\int_{B_{\si(f)}(\b(f))}f\mrm dV_g=\tau\ge\d\q\q\q\int_{\Si\sm B_{\si(f)}(\b(f))}f\mrm dV_g=1-\tau\ge\d.$$
\eite
\eite\

To prove the final assertion, write (up to sub-sequences), $(\b_\infty,\si_\infty)=\lim_{n\to+\infty}(\b(f^n),\si(f^n))$.\\
For large $n$ we will have
$$\int_{\Si\sm\Cup_{x'\in\mcal F}B_\d(x')}f^n\mrm dV_g\le\fr{\d}2\q\q\q\int_{B_\d(x'')}f^n\mrm dV_g\le\fr{\d}2\q\q\q\tx{for any }x''\in\mcal F\sm\{x\},$$
which excludes $\si_\infty=\d$. We also exclude $\si_\infty\in(0,\d)$ as it would give
$$\int_{B_{\fr{3}2\si_\infty}(\b_\infty)}f^n\mrm dV_g\ge\d\q\q\q\int_{\Si\sm B_\fr{\si_\infty}2(\b_\infty)}f^n\mrm dV_g\ge\d.$$
which is a contradiction since $\mcal F\cap\l(\ol{A_{\fr{\si_\infty}2,\fr{3}2\si_\infty}(\b_\infty)}\r)=\es$.\\
Finally, we exclude $\b_\infty\ne x$ because we would get the following contradiction:
$$\int_{B_\d(\b_\infty)}f^n\mrm dV_g\ge\d.$$

The number $\tau$ in the proof of Lemma \ref{center} will be chosen later in Section \ref{s:ex} in such a way that it verifies some good properties when evaluated on the test functions constructed in Section \ref{s:test}.
\epf

Combining such a map $\psi$ with Lemma \ref{mtmacro} we deduce some extra information on low sub-levels of $J_\rho$.\\

\bcor\label{scale}
Let $\d,\psi$ be as in Lemma \ref{center} and define, for $u\in H^1(\Si)^2$,
$$\b_1(u)=\b_{\Si_{\rho_1,\ul\a_1}}(f_{1,u}),\q\q\si_1(u)=\si_{\Si_{\rho_1,\ul\a_1}}(f_{2,u})\q\q\q\b_2(u)=\b_{\Si_{\rho_2,\ul\a_2}}(f_{2,u}),\q\q\si_2(u)=\si_{\Si_{\rho_2,\ul\a_2}}(f_{2,u}).$$
Then for any $\d'>0$ there exists $L_{\d'}$ such that if $\si_i(u)\ge\d'$ for both $i=1,2$, then $J_\rho(u)\ge-L_{\d'}$.
\ecor

\bpf
Assume first $\si_1(u)=\d$: from the statement of Lemma \ref{scale}, we get one of the following:
\bite
\item $\int_{\Si\sm\Cup_{p=1}^MB_\d(p_m)}f_{1,u}\mrm dV_g\ge\fr{\d}2$,
\item $\int_{B_\d(p_m)}f_{1,u}\mrm dV_g\ge\fr{\d}{2M}$ for some $p_m\nin\Si_{\rho_1,\ul\a_1}$,
\item $\int_{B_\d(p_m')}f_{1,u}\mrm dV_g\ge\d,\,\int_{B_\d(p_{m''})}f_{1,u}\mrm dV_g\ge\d$ for some $m'\ne m''$.
\eite
Depending on which possibility occurs, define respectively
\bite
\item $\O_{11}:=\Si\sm\Cup_{p=1}^MB_\d(p_m)$,
\item $\O_{11}:=B_\d(p_m)$,
\item $\O_{11}:=B_\d(p_{m'}),\,\O_{12}:=B_\d(p_{m''})$.
\eite
It is easy to verify that such sets satisfy the hypotheses of Lemma \ref{mtmacro}, up to eventually redefining the map $\psi$ with a smaller $\d\le\fr{\min_{m\ne m'}d(p_m,p_{m'})}4$: in the first case, we have $J_1=0,K_1=1$, in the second case either $J_1=0,K_1=1$ or $J_1=1,K_1=0$ but $\rho<4\pi(1+\a_{1m})$, and in the third case we have $J_1=2,K_1=0$.\\
If $\d'\le\si_1(u)<\d$, then $\int_{\Si\sm B_{\d'}(\b_1(u))}f_{1,u}\mrm dV_g\ge\d$, so we have one between the following:
\bite
\item $\int_{\Si\sm\Cup_{m=1}^MB_\d(x)}f_{1,u}\mrm dV_g\ge\fr{\d}2$
\item $\int_{B_\d(\b_1(u))}f_{1,u}\mrm dV_g\ge\d,\,\int_{B_\d(p_m)}f_{1,u}\mrm dV_g\ge\fr{\d}{2M}$ for some $p_m\ne\b_1(u)$.
\item $\int_{A_{\d',\d}(\b_1(u))}f_{1,u}\mrm dV_g$.
\eite
Depending on which is the case, define:
\bite
\item $\O_{11}:=\Si\sm\Cup_{m=1}^MB_\d(p_m).$
\item $\O_{11}:=B_\d(u)(\b_1(u)),\,\O_{12}:=B_\d(p_m)$.
\item $\O_{11}:=A_{\d',\d}(\b_1(u))$
\eite
Repeat the same argument for $u_2$ to get similarly $\O_{21}$, and possibly $\O_{22}$. Now apply Lemma \ref{mtmacro} and you will get $J_\rho(u)\ge-L_{\d'}$.
\epf

In Section \ref{s:impr}, we will need to combine different types of improved Moser-Trudinger inequalities. To do this, we will need the following technical estimates concerning averages of functions on balls and their boundary:\\

\blem\label{average}
There exists $C>0$ such that for any $u\in H^1(\Si),\,x\in\Si,\,r>0$ one has
$$\l|\fint_{B_r(x)}u\mrm dV_g-\fint_{\pa B_r(x)}u\mrm dV_g\r|\le C\sqrt{\int_{B_r(x)}|\n u|^2\mrm dV_g}.$$
Moreover, for any $R>1$ there exists $C=C_R$ such that
$$\l|\fint_{B_r(x)}u\mrm dV_g-\fint_{B_{Rr}(x)}u\mrm dV_g\r|\le C\sqrt{\int_{B_r(x)}|\n u|^2\mrm dV_g}.$$
The same inequalities hold if $B_r(x)$ is replaced by a domain $\O\sub B_{Rr}(x)$ such that $\O\in\mfrak A_{\d r}$ for some $\d>0$, with $C$ and $C_R$ replaced by some $C_\d,C_{R,\d}>0$, respectively.\\

\elem

The proof of the above lemma follows from the Poincaré-Wirtinger and trace inequalities, which are invariant by dilation. Details can be found, for instance, in \cite{gt}. We will also need the following estimate on harmonic liftings.\\

\blem\label{trace}
Let $r_2>r_1>0$, $f\in H^1(B_{r_2}(0))$ with $\int_{B_{r_1}(0)}f(x)dx=0$ be given and $u$ be the solution of
$$\l\{\bll-\D u=0&\tx{in } A_{r_1,r_2}(0)\\
u=f&\tx{on }\pa B_{r_1}(0)\\u=0&\tx{on }\pa B_{r_2}(0)\earr\r..$$
Then, there exists $C=C_\fr{r_2}{r_1}>0$ such that
$$\int_{A_{r_1,r_2}(0)}|\n u(x)|^2\mrm dx\le C\int_{A_{r_1,r_2}(0)}|\n f(x)|^2\mrm dx$$
\elem

Again, the proof uses elementary techniques in elliptic PDEs, such as Dirichlet principle and Poincar\'e inequality, hence can be found in most textbooks.

\section{The topology of the space $\mcal X$}\label{s:top}

Let us introduce the space $\mcal X$, which will play a fundamental role in all the rest of the paper.
It is obtained removing some points from the join of the weighted barycenters $\Si_{\rho_1,\ul\a_1}\star\Si_{\rho_2,\ul\a_2}$ defined by \eqref{wbar} and \eqref{join}. The points to exclude correspond to improved inequalities for functions centered around the same point and at the same rate of concentration (see Section \ref{s:impr} for more details).\\
Precisely, we have:
\bequ\label{x}
\mcal X:=\Si_{\rho_1,\ul\a_1}\star\Si_{\rho_2,\ul\a_2}\sm\l\{\l(p_m,p_m,\fr{1}2\r):\,\rho_1,\rho_2<4\pi(2+\a_{1m}+\a_{2m})\r\}.
\eequ
In this section, we will prove that, under the assumptions of Theorem \ref{ex}, the space $\mcal X$ is not contractible. In particular, we will prove that it has a non-trivial homology group.\\

In order to do this, we will recall how to calculate the homology groups of the join of two known spaces. Since the join is homotopically equivalent to a smash product of $X,\,Y$ and $\S^1$ (see \cite{hat} for details), its homology groups only depend on the homology of $X$ and $Y$.\\

\bthm\label{hom}(\cite{hat}, Theorem $3.21$)
Let $X$ and $Y$ be two topological spaces. Then,
$$\wt H_q(X\star Y)=\sum_{q'=0}^q\wt H_{q'}(X)\pl\wt H_{q-q'-1}(Y).$$
In particular, if $X=\l(\S^{D_1}\r)^{\vee N_1}$ and $Y=\l(\S^{D_2}\r)^{\vee N_2}$ are wedge sum of spheres, then $X\star Y$ has the same homology of $\l(\S^{D_1+D_2+1}\r)^{\vee N_1N_2}$.
\ethm

Actually, in the same book \cite{hat} it is shown that the following homotopical equivalence holds: $\l(\S^{D_1}\r)^{\vee N_1}\star\l(\S^{D_2}\r)^{\vee N_2}\simeq\l(\S^{D_1+D_2+1}\r)^{\vee N_1N_2}$.\\

Here is the main result of this section:\\

\bthm\label{top}
Let $M_1,M_2,M_3$ be as in \eqref{m123} and $\mcal X$ be as in \eqref{x} and suppose
\bequ\label{m1m2m3}
(M_1,M_2,M_3)\nin\{(1,m,0),(m,1,0),(2,2,1),(2,3,2),(3,2,2),\,m\in\N\}.
\eequ
Then, the space $\mcal X$ has non-trivial homology groups. In particular, it is not contractible.
\ethm

The assumptions on the $M_1,M_2,M_3$, that is, respectively on the cardinality of $\Si_{\rho_1,\ul\a_1},\,\Si_{\rho_2,\ul\a_2}$ and on the number of midpoints to be removed, are actually sharp.\\
This can be seen clearly from the Figure $1$: the configurations $M_1=1,\,M_3=0$ are star-shaped, and even in the two remaining case it is easy to see $\mcal X$ has trivial topology. On the other hand, Figure $2$ shows a non-contractible configuration.\\

\begin{figure}[h!]
\center

\includegraphics[width=0.2\linewidth]{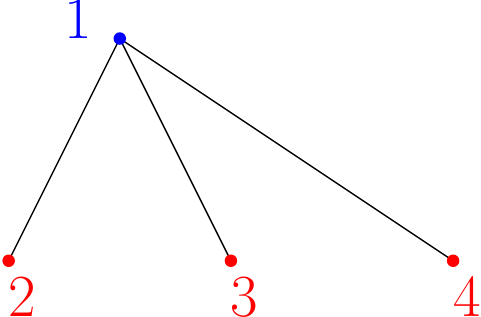}\q\q\q\q\q\q
\includegraphics[width=0.1\linewidth]{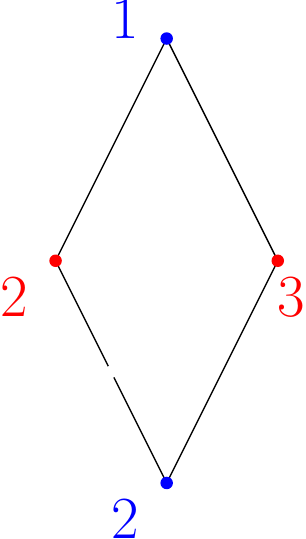}\q\q\q\q\q\q
\includegraphics[width=0.2\linewidth]{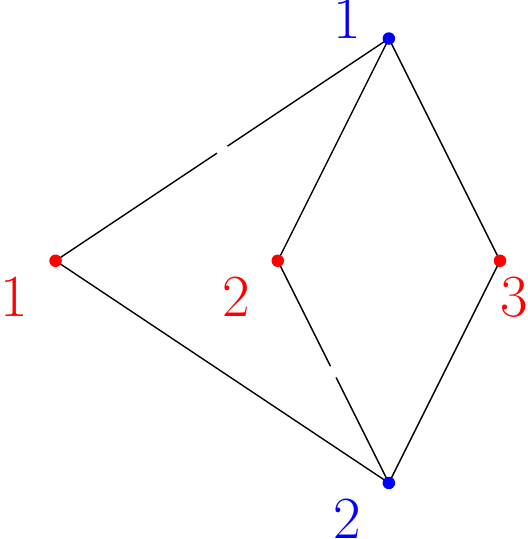}\\

\caption{The space $\mcal X$ in the cases $(M_1,M_2,M_3)\in\{(1,3,0),(2,2,1),(2,3,2)\}$ (contractible).}\

\includegraphics[width=0.3\linewidth]{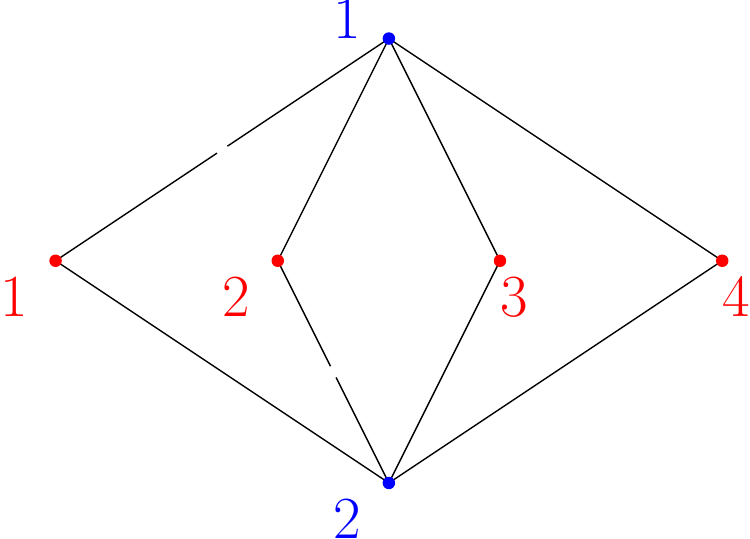}\\
\caption{The space $\mcal X$ in the case $M_1=2,\,M_2=4,\,M_3=2$ (not contractible).}
\end{figure}

\bpf[Proof of Theorem \ref{top}]
The spaces $\Si_{\rho_i,\ul\a_i}$ are discrete sets of $M_i$ points, for $i=1,2$, that is a wedge sum of $M_i-1$ copies of $\S^0$. Therefore, by Theorem \ref{hom}, $\Si_{\rho_1,\ul\a_1}\star\Si_{\rho_2,\ul\a_2}$ has the same homology as $\l(\S^1\r)^{\vee(M_1-1)(M_2-1)}$.\\
The set we have to remove from the join is made up by $M_3$ singular points $\{p_{m_1},\ds,p_{m_{M_3}}\}$ for some $\{m_1,\ds,m_{M_3}\}\sub\{1,\ds,M\}$.\\
Defining then, for some fixed $\d<\fr{1}2$, $\mcal Y:=\Cup_{j=1}^{M_3}B_\d\l(p_{m_j},p_{m_j},\fr{1}2\r)$,
$\mcal Y$ retracts on $\l\{p_{m_1},\ds,p_{m_{M_3}}\r\}$. On the other hand, $\mcal X\cap\mcal Y$ is a disjoint union of $M_3$ punctured intervals, that is a discrete set of $2M_3$ points, and $\mcal X\cup\mcal Y$ is the whole join.
Therefore, the Mayer-Vietoris sequence yields
$$\ub{H_1(\mcal X\cap\mcal Y)}_0\to H_1(\mcal X)\pl\ub{H_1(\mcal Y)}_0\to\ub{H_1(\mcal X\cup\mcal Y)}_{\Z^{(M_1-1)(M_2-1)}}\to\ub{\wt H_0(\mcal X\cap\mcal Y)}_{\Z^{2M_3-1}}\to\wt H_0(\mcal X)\pl\ub{\wt H_0(\mcal Y)}_{\Z^{M_3-1}}\to\ub{\wt H_0(\mcal X\cup\mcal Y)}_0.$$
The exactness of the sequence implies that $b_1(\mcal X)-\wt b_0(\mcal X)=(M_1-1)(M_2-1)-M_3$, so if the latter number is not zero we get at least a non-trivial homology group.\\
Algebraic computations show that, under the assumption $M_1,M_2\ge M_3$, $(M_1-1)(M_2-1)\ne M_3$ is equivalent to \eqref{m1m2m3}, therefore the proof is complete.
\epf

\section{Construction of test functions}\label{s:test}

We will now introduce some test functions from the space $\mcal X$, introduced in Section \ref{s:top}, to arbitrarily low sub-levels. Such test functions will have a profile which resembles the entire solutions of the Liouville equation and of the Toda system: it will not always suffice to consider the \emph{standard} bubbles
$$\ph^\la_{p,\a}=-2\log\max\l\{1,(\la d(\cd,p))^{2(1+\a)}\r\},$$
which roughly resemble the solutions of the scalar Liouville equation. This is because, when the two components are centered at the same points, a higher amount of energy is needed due to the expression of $Q$ (see the Introduction) which penalizes parallel gradients. . This is basically the reason that the join is punctured in \eqref{x}.\\
We will need two more profiles for the construction of $\Phi^\la$, which have been considered in \cite{jkm} for the regular Toda system.
$${\ph'}^\la_{p,\a_1,\a_2}=-2\log\max\l\{1,(\la d(\cd,p))^{2(2+\a_1+\a_2)}\r\},$$
$${\ph''}^\la_{p,\a_1,\a_2}=-2\log\max\l\{1,\la^{2(2+\a_1+\a_2)}d(\cd,p)^{2(1+\a_1)}\r\}.$$
We will use suitable interpolation between each of the above three profiles depending on whether the points $x_i\in\Si_{\rho_i,\ul\a_i}$ coincide or not and depending on which of the parameters $\rho_i$ is greater or less than $4\pi\l(\o_{\ul\a_1}(\{x_i\})+\o_{\ul\a_2}(\{x_i\})\r)$, see \eqref{wbar}. The map $\Phi^\la$ will therefore be defined case by case, hence its definition will be quite lengthy and will be postponed in the proof of the theorem, rather than in its statement.\\
As a final remark, we considered a \emph{truncated} version of the bubbles instead of the usual smooth ones. Under this change we get very similar estimates, though with simpler calculations, since truncated functions are easy to handle.\\

\bthm\label{test}
There exists a family of maps $\l\{\Phi^\la\r\}_{\la>2}:\mcal X\to H^1(\Si)^2$ such that
$$J_\rho\l(\Phi^\la(\z)\r)\us{\la\to+\infty}\lto-\infty\q\q\q\tx{uniformly for }\z\in\mcal X.$$
\ethm

\bpf

Let us start by defining $\Phi^\la(\z)=\l(\ph_1-\fr{\ph_2}2,\ph_2-\fr{\ph_1}2\r)$ when $\z=(p_m,p_m,t)$ for some $m$. $\Phi^\la$ will be defined in different ways, depending on the relative positions of $\rho_1,\rho_2,\a_{1m},\a_{2m}$ in $\R$.\\

\bite
\item[$(<<)$] $\rho_1,\rho_2<4\pi(2+\a_{1m}+\a_{2m})$:
$$\bl
\ph_1:=\l\{\bll-2\log\max\l\{1,(\la d(\cd,p_m))^{2(1+\a_{1m})}\r\}&\tx{if }t<\fr{1}2\\0&\tx{if }t>\fr{1}2\earr\r.\\
\ph_2:=\l\{\bll0&\tx{if }t<\fr{1}2\\-2\log\max\l\{1,(\la d(\cd,p_m))^{2(1+\a_{2m})}\r\}&\tx{if }t>\fr{1}2.\earr\r.\earr$$

\item[$(<>)$] $\rho_1<4\pi(2+\a_{1m}+\a_{2m})<\rho_2$:
$$\bl\ph_1:=-2\log\max\l\{1,\max\l\{1,(\la t)^{2(1+\a_{2m})}\r\}(\la d(\cd,p_m))^{2(1+\a_{1m})}\r\}\\\ph_2:=-2\log\max\l\{1,(\la td(\cd,p_m))^{2(2+\a_{1m}+\a_{2m})}\r\}.\earr$$

\item[$(><)$] $\rho_2<4\pi(2+\a_{1m}+\a_{2m})<\rho_1$:
$$\bl\ph_1:=-2\log\max\l\{1,(\la(1-t)d(\cd,p_m))^{2(2+\a_{1m}+\a_{2m})}\r\}\\\ph_2:=-2\log\max\l\{1,\max\l\{1,(\la(1-t))^{2(1+\a_{1m})}\r\}(\la d(\cd,p_m))^{2(1+\a_{2m})}\r\}.\earr$$

\item[$(>>)$] $\rho_1,\rho_2>4\pi(2+\a_{1m}+\a_{2m})$:
$$\bl\ph_1:=-2\log\max\l\{1,\l(\la\fr{\max\{1,\la t\}}{\max\{1,\la(1-t)\}}\r)^{2+\a_{1m}+\a_{2m}}d(\cd,p_m)^{2(1+\a_{1m})},(\la d(\cd,p_m))^{2(2+\a_{1m}+\a_{2m})}\r\}\\
\ph_2:=-2\log\max\l\{1,\l(\la\fr{\max\{1,\la(1-t)\}}{\max\{1,\la t\}}\r)^{2+\a_{1m}+\a_{2m}}d(\cd,p_m)^{2(1+\a_{2m})},(\la d(\cd,p_m))^{2(2+\a_{1m}+\a_{2m})}\r\}.\earr$$

\eite
We will need some estimates on $\Phi^\la$, which will be proved in three separates lemmas and which, combined, will give the proof of the theorem.
\epf

\tbf{Convention:} When using normal coordinates near the peaks of the test functions, the metric coefficients will slightly deviate from the Euclidean ones. We will then have coefficients of order $(1+o_\la(1))$ in front of the logarithmic terms appearing below. To keep the formulas shorter, we will omit them, as they will be harmless for the final estimates.\\

\blem\label{grad}
Let $\ph_1,\ph_2$ be as in Theorem \ref{test}. Then, setting $\mcal Q:=\int_\Si Q\l(\ph_1-\fr{\ph_2}2,\ph_2-\fr{\ph_1}2\r)\mrm dV_g$,
in each case we have for $\la$ large

\bite
\item[$(<<)$]\q$\mcal Q=\l\{\bll8\pi(1+\a_{1m})^2\log\la+O(1)&\tx{if }t<\fr{1}2\\8\pi(1+\a_{2m})^2\log\la+O(1)&\tx{if }t>\fr{1}2\earr\r.$.

\item[$(<>)$]\q$\mcal Q=8\pi(2+\a_{1m}+\a_{2m})^2\log\max\{1,\la t\}+8\pi(1+\a_{1m})^2\log\min\l\{\la,\fr{1}t\r\}+O(1)$.

\item[$(><)$]\q$\mcal Q=8\pi(2+\a_{1m}+\a_{2m})^2\log\max\{1,\la(1-t)\}+8\pi(1+\a_{2m})^2\log\min\l\{\la,\fr{1}{1-t}\r\}+O(1)$.

\item[$(>>)$]\q$\mcal Q=8\pi(2+\a_{1m}+\a_{2m})^2\log\la+O(1)$.

\eite
\elem

\bpf
Let us start by the case $(<<)$. We assume $t<\fr{1}2$, since the case $t>\fr{1}2$ can be treated in the very same way just switching the indices.
There holds
$$\n\ph_1=\l\{\bll0&\tx{if }d(\cd,p_m)<\fr{1}\la\\-4(1+\a_{1m})\fr{\n d(\cd,p_m)}{d(\cd,p_m)}&\tx{if }d(x,p_m)>\fr{1}\la\earr\r.\q\q\q\n\ph_2\eq0.$$
Therefore, since $|\n d(\cd,p_m)|=1$ a.e. on $\Si$, we get
\bequ
\mcal Q=\fr{1}4\int_\Si|\n\ph_1|^2\mrm dV_g=4(1+\a_{1m})^2\int_{\Si\sm B_\fr{1}\la(p_m)}\fr{\mrm dV_g}{d(\cd,p_m)^2}=8\pi(1+\a_{1m})^2\log\la+O(1).
\eequ
In the case $(<>)$ we can assume $\la t\ge1$, since otherwise $\ph_1,\ph_2$ are defined just like the previous section.
We have
$$\n\ph_1=\l\{\bll0&\tx{if }d(\cd,p_m)<\fr{1}{\la(\la t)^\fr{1+\a_{2m}}{1+\a_{1m}}}\\-4(1+\a_{1m})\fr{\n d(\cd,p_m)}{d(\cd,p_m)}&\tx{if }d(x,p_m)>\fr{1}{\la(\la t)^\fr{1+\a_{2m}}{1+\a_{1m}}}\earr\r.,$$
$$\n\ph_2=\l\{\bll0&\tx{if }d(\cd,p_m)<\fr{1}{\la t}\\-4(2+\a_{1m}+\a_{2m})\fr{\n d(\cd,p_m)}{d(\cd,p_m)}&\tx{if }d(x,p_m)>\fr{1}{\la t}\earr\r.$$
therefore
\beqy
\mcal Q&=&\fr{1}4\int_\Si|\n\ph_1|^2\mrm dV_g-\fr{1}4\int_\Si\n\ph_1\cd\n\ph_2\mrm dV_g+\fr{1}4\int_\Si|\n\ph_2|^2\mrm dV_g\\
&=&4(1+\a_{1m})^2\int_{\Si\sm B_\fr{1}{\la(\la t)^\fr{1+\a_{2m}}{1+\a_{1m}}}(p_m)}\fr{\mrm dV_g}{d(\cd,p_m)^2}-4(2+\a_{1m}+\a_{2m})(1+\a_{1m})\int_{\Si\sm B_\fr{1}{\la t}(p_m)}\fr{\mrm dV_g}{d(\cd,p_m)^2}\\
&+&4(2+\a_{1m}+\a_{2m})^2\int_{\Si\sm B_\fr{1}{\la t}(p_m)}\fr{\mrm dV_g}{d(\cd,p_m)^2}\\
&=&8\pi(1+\a_{1m})^2\log\la+8\pi(1+\a_{1m})(1+\a_{2m})\log(\la t)\\
&-&8\pi(2+\a_{1m}+\a_{2m})(1+\a_{1m})\log(\la t)+8\pi(2+\a_{1m}+\a_{2m})^2\log(\la t)+O(1)\\
&=&8\pi(2+\a_{1m}+\a_{2m})^2\log(\la t)+4\pi(1+\a_{2m})^2\log\fr{1}t+O(1).
\eeqy
In the case $(><)$ we can argue as in $(<>)$ just switching the indices; similar calculations also yield the last case $(>>)$.
\epf

\blem\label{averphi}
Let $\ph_1,\ph_2$ be as above. Then, in each case we have:\\
\bite
\item[$(<<)$]
$$\int_\Si\ph_1\mrm dV_g=\l\{\bll-4(1+\a_{1m})\log\la+O(1)&\tx{if }t<\fr{1}2\\0&\tx{if }t>\fr{1}2\earr\r.,\q\int_\Si\ph_2\mrm dV_g=\l\{\bll0&\tx{if }t<\fr{1}2\\-4(1+\a_{2m})\log\la+O(1)&\tx{if }t>\fr{1}2.\earr\r.$$
\item[$(<>)$]
$$\int_\Si\ph_1\mrm dV_g=-4(1+\a_{1m})\log\la-4(1+\a_{2m})\log\max\{1,\la t\}+O(1),$$
$$\int_\Si\ph_2\mrm dV_g=-4(2+\a_{1m}+\a_{2m})\log\max\{1,\la t\}+O(1).$$
\item[$(><)$]
$$\int_\Si\ph_1\mrm dV_g=-4(2+\a_{1m}+\a_{2m})\log\max\{1,\la(1-t)\}+O(1),$$
$$\int_\Si\ph_2\mrm dV_g=-4(1+\a_{2m})\log\la-4(1+\a_{1m})\log\max\{1,\la(1-t)\}+O(1).$$
\item[$(>>)$]
$$\int_\Si\ph_1\mrm dV_g=\int_\Si\ph_2\mrm dV_g+O(1)=-4(2+\a_{1m}+\a_{2m})\log\la+O(1).$$
\eite
\elem

\bpf
Let us consider the case $(<<),\,t<\fr{1}2$. Since we have
$$-4(1+\a_{1m})\log\max\{1,d(\cd,p_m)\}\le\ph_1+4(1+\a_{1m})\log\la\le-4(1+\a_{1m})\log d(\cd,p_m),$$
with both the first and the last function having finite average over $\Si$, we are done.\\

The same argument also works in all the other cases.
\epf

\blem\label{logint}
Let $\ph_1,\ph_2$ be as above. Then, in each case we have:\\

\bite
\item[$(<<)$]
$$\log\int_\Si\wt h_1e^{\ph_1-\fr{\ph_2}2}\mrm dV_g=\l\{\bll-2(1+\a_{1m})\log\la+O(1)&\tx{if }t<\fr{1}2\\2(1+\a_{2m})\log\la+O(1)&\tx{if }t>\fr{1}2\earr\r.,$$
$$\log\int_\Si\wt h_2e^{\ph_2-\fr{\ph_1}2}\mrm dV_g=\l\{\bll2(1+\a_{1m})\log\la+O(1)&\tx{if }t<\fr{1}2\\-2(1+\a_{2m})\log\la+O(1)&\tx{if }t>\fr{1}2.\earr\r..$$
\item[$(<>)$]
$$\log\int_\Si\wt h_1e^{\ph_1-\fr{\ph_2}2}\mrm dV_g=-2(1+\a_{1m})\log\la-2(1+\a_{2m})\log\max\{1,\la t\},$$
$$\log\int_\Si\wt h_2e^{\ph_2-\fr{\ph_1}2}\mrm dV_g=2(1+\a_{1m})\min\l\{\la,\fr{1}t\r\}.$$
\item[$(><)$]
$$\log\int_\Si\wt h_1e^{\ph_1-\fr{\ph_2}2}\mrm dV_g=2(1+\a_{2m})\min\l\{\la,\fr{1}{1-t}\r\},$$
$$\log\int_\Si\wt h_2e^{\ph_2-\fr{\ph_1}2}\mrm dV_g=-2(1+\a_{2m})\log\la-2(1+\a_{1m})\log\max\{1,\la(1-t)\}.$$
\item[$(>>)$]
$$\log\int_\Si\wt h_1e^{\ph_1-\fr{\ph_2}2}\mrm dV_g=-(2+\a_{1m}+\a_{2m})\log\l(\la\fr{\max\{1,\la t\}}{\max\{1,\la(1-t)\}}\r)+O(1)$$
$$\log\int_\Si\wt h_2e^{\ph_2-\fr{\ph_1}2}\mrm dV_g=-(2+\a_{1m}+\a_{2m})\log\l(\la\fr{\max\{1,\la(1-t)\}}{\max\{1,\la t\}}\r)+O(1)$$
\eite
\elem

\bpf
Again, we will just consider the first case.\\
Given any $\d>0$, if $d(\cd,p_m)\ge\d$ one has
$$e^{\ph_1-\fr{\ph_2}2}\ge\fr{C_\d}{\la^{4(1+\a_{1m})}}\q\q\q e^{\ph_1-\fr{\ph_2}2}\ge C_\d\la^{2(1+\a_{1m})},$$
therefore we will suffice to consider only the integral on $B_\d(p_m)$:
\beqy
\log\int_{B_\d(p_m)}\wt h_1e^{\ph_1-\fr{\ph_2}2}\mrm dV_g&=&\log\int_{B_\d(p_m)}d(\cd,p_m)^{2\a_{1m}}e^{\ph_1-\fr{\ph_2}2}\mrm dV_g+O(1)\\
&=&\log\l(\ub{\int_{B_\fr{1}\la(p_m)}d(\cd,p_m)^{2\a_{1m}}\mrm dV_g}_{\sim\la^{-2(1+\a_{1m})}}+\la^{-4(1+\a_{1m})}\ub{\int_{A_{\fr{1}\la,\d}(p_m)}\fr{\mrm dV_g}{d(\cd,p_m)^{2(2+\a_{1m})}}}_{\sim\la^{2(1+\a_{1m})}}\r)+O(1)\\
&=&-2(1+\a_{1m})\log\la+O(1),
\eeqy
\beqy
\log\int_{B_\d(p_m)}\wt h_1e^{\ph_2-\fr{\ph_1}2}\mrm dV_g&=&\log\l(\ub{\int_{B_\fr{1}\la(p_m)}d(\cd,p_m)^{2\a_{2m}}\mrm dV_g}_{\sim\la^{-2(1+\a_{2m})}}+\la^{2(1+\a_{1m})}\ub{\int_{A_{\fr{1}\la,\d}(p_m)}d(\cd,p_m)^{2(1+\a_{1m}+\a_{2m})}\mrm dV_g}_{\sim1}\r)+O(1)\\
&=&2(1+\a_{1m})\log\la+O(1),
\eeqy
\epf

\bpf[Proof of Theorem \ref{test}, continued]
From the previous lemmas we can easily prove the theorem in the case $x_1=x_2$. In fact, writing
$$J_\rho\l(\ph_1-\fr{\ph_2}2,\ph_2-\fr{\ph_1}2\r)=\int_\Si Q\l(\ph_1-\fr{\ph_2}2,\ph_2-\fr{\ph_1}2\r)\mrm dV_g-\sum_{i=1}^2\rho_i\l(\log\int_\Si\wt h_ie^{\ph_i-\fr{\ph_{3-i}}2}\mrm dV_g-\int_\Si\ph_i\mrm dV_g+\fr{1}2\int_\Si\ph_{3-i}\mrm dV_g\r),$$
we get, in each case,
\bite
\item[$(<<)$]
$$J_\rho\l(\ph_1-\fr{\ph_2}2,\ph_2-\fr{\ph_1}2\r)=\l\{\bll2(1+\a_{1m})(4\pi(2+\a_{1m})-\rho_1)\log\la+O(1)&\tx{if }t<\fr{1}2\\2(1+\a_{2m})(4\pi(2+\a_{2m})-\rho_2)\log\la+O(1)&\tx{if }t>\fr{1}2\earr\r.,$$
\item[$(<>)$]\beqy
J_\rho\l(\ph_1-\fr{\ph_2}2,\ph_2-\fr{\ph_1}2\r)&=&2(1+\a_{1m})(4\pi(1+\a_{1m})-\rho_1)\log\min\l\{\la,\fr{1}t\r\}\\
&+&2(2+\a_{1m}+\a_{2m})(4\pi(2+\a_{1m}+\a_{2m})-\rho_2)\log\max\{1,\la t\}+O(1),
\eeqy
\item[$(><)$]
\beqy
J_\rho\l(\ph_1-\fr{\ph_2}2,\ph_2-\fr{\ph_1}2\r)&=&2(2+\a_{1m}+\a_{2m})(4\pi(2+\a_{1m}+\a_{2m})-\rho_1)\log\max\{1,\la(1-t)\}\\
&+&2(1+\a_{2m})(4\pi(1+\a_{2m})-\rho_2)\log\min\l\{\la,\fr{1}{1-t}\r\}+O(1),
\eeqy
\item[$(>>)$]
\beqy J_\rho\l(\ph_1-\fr{\ph_2}2,\ph_2-\fr{\ph_1}2\r)&=&(2+\a_{1m}+\a_{2m})(4\pi(2+\a_{1m}+\a_{2m})-\rho_1)\log\l(\la\fr{\max\{1,\la(1-t)\}}{\max\{1,\la t\}}\r)\\
&+&(2+\a_{1m}+\a_{2m})(4\pi(2+\a_{1m}+\a_{2m})-\rho_2)\log\l(\la\fr{\max\{1,\la t\}}{\max\{1,\la(1-t)\}}\r)+O(1),
\eeqy
\eite
which all tend to $-\infty$ independently of $t$.\\

Let us now consider the case $x_1\ne x_2$.\\
Here, we define $\Phi^\la$ just by interpolating linearly between the test functions defined before:
$$\Phi^\la(x_1,x_2,t)=\Phi^{\la(1-t)}(x_1,x_1,0)+\Phi^{\la t}(x_2,x_2,1).$$
Since $d(p_m,p_{m'})\ge\d>0$, then the bubbles centered at $p_m$ and $p_{m'}$ do not interact, therefore the estimates from Lemmas \ref{grad}, \ref{averphi}, \ref{logint} also work for such test functions. We will show this fact in detail in the case $\rho_1,\rho_2<4\pi(2+\a_{1m}+\a_{2m}),4\pi(2+\a_{1m'}+\a_{2m'})$. Writing
$$(\ph_1,\ph_2)=\l(-2\log\max\l\{1,(\la(1-t)d(\cd,p_m))^{2(1+\a_{1m})}\r\},-2\log\max\l\{1,(\la td(\cd,p_{m'}))^{2(1+\a_{2m'})}\r\}\r),$$
by the previous explicit computation of $\n\ph_1,\n\ph_2$ we get
\beqa
\nonumber\mcal Q&=&\fr{1}4\int_{B_\d(p_m)}|\n\ph_1|^2\mrm dV_g+\fr{1}4\int_{B_\d(p_{m'})}|\n\ph_2|^2\mrm dV_g+O(1)\\
\label{phi1}&=&8\pi(1+\a_{1m})^2\log\max\{1,\la(1-t)\}+8\pi(2+\a_{2m'})^2\log\max\{1,\la t\}+O(1).
\eeqa
Moreover, by linearity,
\bequ\label{phi2}
\int_\Si\ph_1\mrm dV_g=-4(1+\a_{1m})\log\max\{1,\la(1-t)\}+O(1)\q\int_\Si\ph_2\mrm dV_g=-4(1+\a_{2m'})\log\max\{1,\la t\}+O(1).
\eequ
Finally, as before the integral of $\wt h_1e^{\ph_1-\fr{\ph_2}2}$ is negligible outside $B_\d(p_m)$, and inside the ball we have $\fr{1}{C_\d}\le\l|\ph_2-\int_\Si\ph_2\mrm dV_g\r|\le C$ on $B_\d(p_m)$, hence
\beqa
\nonumber\log\int_\Si\wt h_1e^{\ph_1-\fr{\ph_2}2}\mrm dV_g&=&\log\l(\max\{1,\la t\}^{2(1+\a_{2m'})}\int_{B_\fr{1}{\max\{1,\la(1-t)\}}(p_m)}d(\cd,p_m)^{2\a_{1m}}\mrm dV_g\r.\\
\nonumber&+&\l.\max\{1,\la t\}^{2(1+\a_{2m'})}\max\{1,\la(1-t)\}^{2(1+\a_{1m})}\int_{A_{\fr{1}\la,\d}(p_m)}\fr{\mrm dV_g}{d(\cd,p_m)^{2(2+\a_{1m})}}\r)+O(1)\\
\label{phi3}&=&2(1+\a_{2m'})\log\max\{1,\la t\}-2(1+\a_{1m})\log\max\{1,\la(1-t)\}+O(1)
\eeqa
and similarly
$$\log\int_\Si\wt h_2e^{\ph_2-\fr{\ph_1}2}\mrm dV_g=2(1+\a_{1m})\log\max\{1,\la(1-t)\}-2(1+\a_{1m'})\log\max\{1,\la t\}+O(1).$$
Therefore, by \eqref{phi1}, \eqref{phi2} and \eqref{phi3} we deduce
\beqy
J_\rho\l(\ph_1-\fr{\ph_2}2,\ph_2-\fr{\ph_1}2\r)&=&2(1+\a_{1m})(4\pi(1+\a_{1m})-\rho_1)\log\max\{1,\la(1-t)\}\\
&+&2(1+\a_{2m'})(4\pi(1+\a_{2m'})-\rho_2)\log\max\{1,\la t\}+O(1).
\eeqy
This concludes the proof.
\epf

\section{Improved Moser-Trudinger inequalities}\label{s:impr}

In this section we will deduce some improved Moser-Trudinger inequalities when the two components have the same center and mass of concentration, in the sense defined by Lemma \ref{center}.\\

\bthm\label{impr}
Let $\beta_i(u),\si_i(u)$ be as in Corollary \ref{scale}. There exists $L\gg0$ such that if
$$\l\{\bll\b_1(u)=\b_2(u)=p_m&\tx{with }\rho_1,\rho_2<4\pi(2+\a_{1m}+\a_{2m})\\\si_1(u)=\si_2(u)&\earr\r.,$$
then $J_\rho(u)\ge-L$.
\ethm

Theorem \ref{impr} is based on the following two lemmas, inspired by \cite{mr13}.
Basically, we assume $u_1$ and $u_2$ to have the same center and scale of concentration and we provide local estimates in a ball which is roughly centered at the center of mass and whose radius is roughly the same as the scale of concentration. Inner estimates use a dilation argument, outer estimates use a Kelvin transform. With respect to the above-cited paper, we also have to consider concentration around the boundary of the ball, hence we will combine those arguments with Theorem \ref{mtdisk} and Lemma \ref{mtloc}.\\

\blem\label{mtin}
For any $\e>0,\,\a_1,\a_2\in(-1,0]$ there exists $C=C_\e$ such that for any $p\in\Si,\,s>0$ small enough and $u\in H^1(\Si)^2$ one has
\beqa
\nonumber&&4\pi\sum_{i=1}^2(1+\a_i)\l(\log\int_{B_{\fr{s}2}(p)}d(\cd,p)^{2\a_i}e^{u_i}\mrm dV_g-\fint_{B_s(p)}u_i\mrm dV_g\r)-8\pi\l((1+\a_1)^2+(1+\a_2)^2\r)\log s\\
\label{in}&\le&(1+\e)\int_{B_s(p)}Q(u)\mrm dV_g+C,\\
\nonumber&&4\pi(1+\a_1)\l(\log\int_{B_{\fr{s}8}(p)}d(\cd,p)^{2\a_1}e^{u_1}\mrm dV_g-\fint_{B_s(p)}u_1\mrm dV_g\r)\\
\nonumber&+&2\pi\min\{1,2+\a_1+\a_2\}\l(\log\int_{A_{\fr{s}4,s}(p)}d(\cd,p)^{2\a_2}e^{u_2}\mrm dV_g-\fint_{B_s(p)}u_2\mrm dV_g\r)\\
\nonumber&-&4\pi\l(2(1+\a_1)^2+\min\{1,2+\a_1+\a_2\}(1+\a_2)\r)\log s\\
\label{inmix}&\le&(1+\e)\int_{B_s(p)}Q(u)\mrm dV_g+C,\\
\nonumber&&2\pi\sum_{i=1}^2\min\{1,2+\a_1+\a_2\}\l(\log\int_{A_{\fr{s}2,s}(p)}d(\cd,p)^{2\a_i}e^{u_i}\mrm dV_g-\fint_{B_s(p)}u_i\mrm dV_g\r)\\
\nonumber&-&4\pi\min\l\{2+\a_1+\a_2,(2+\a_1+\a_2)^2\r\}\log s\\
\label{inbdry}&\le&(1+\e)\int_{B_s(p)}Q(u)\mrm dV_g+C.
\eeqa
The last statement holds true if $B_s(p)$ is replaced by $\O_s$ simply connected belonging to $\mfrak A_{\d s}$ (see \eqref{eq:ud}) and such that $B_{\l(\fr{1}2+\d\r)s}(p)\sub\O_s\sub B_{\fr{s}\d}(p)$ for some $\d>0$, with $C$ replaced with some $C_\d>0$.
\elem

\bpf
By assuming $s$ small enough, we can suppose the metric to be flat on $B_s(p)$, up to negligible remainder terms. Therefore, we will assume to work on a Euclidean ball centered at the origin: we will indicate such balls simply as $B_s$, omitting their center, and we will use a similar convention for annuli. Moreover, we will write $|x|$ for $d(x,p)$.\\
Consider the dilation $v_i(z)=u_i(sz)$ for $z\in B_1$. It verifies, for $r\in\l\{\fr{1}8,\fr{1}4,\fr{1}2\r\}$
$$\int_{B_r}|z|^{2\a_i}e^{v_i(z)}\mrm dz=s^{-2-2\a_i}\int_{B_{rs}}|x|^{2\a_i}e^{u_i(x)}\mrm dx,$$
$$\int_{A_{r,1}}|z|^{2\a_i}e^{v_i(z)}\mrm dz=s^{-2-2\a_i}\int_{A_{rs,s}}|x|^{2\a_i}e^{u_i(x)}\mrm dx.$$
$$\int_{B_1}Q(v(z))\mrm dz=\int_{B_s}Q(u(x))\mrm dx,\q\q\q\fint_{B_1}v(z)\mrm dz=\fint_{B_s}u(x)\mrm dx,$$
To get \eqref{in}, it suffices to apply \eqref{loc} to $v=(v_1,v_2)$:
\beqy
&&4\pi\sum_{i=1}^2(1+\a_i)\l(\log\int_{B_{\fr{s}2}}|x|^{2\a_i}e^{u_i(x)}\mrm dx-\fint_{B_s}u_i(x)\mrm dx-2(1+\a_i)\log s\r)\\
&\le&4\pi\sum_{i=1}^2(1+\a_i)\l(\log\int_{B_{\fr{1}2}}|z|^{2\a_i}e^{v_i(z)}\mrm dz-\fint_{B_1}v_i(z)\mrm dz\r)\\
&\le&(1+\e)\int_{B_1}Q(v(z))\mrm dz+C\\
&=&(1+\e)\int_{B_s}Q(u(x))\mrm dx+C.
\eeqy
For \eqref{inmix}, one has to use \eqref{locmix} on $v$, and the elementary fact that $\fr{1}C\le|z|^{2\a_2}\le C$ on $A_{\fr{1}4,1}$:
\beqy
&&4\pi(1+\a_1)\l(\log\int_{B_{\fr{s}8}}|x|^{2\a_1}e^{u_1(x)}\mrm dx-\fint_{B_s}u_1(x)\mrm dx-2(1+\a_1)\log s\r)\\
&+&2\pi\min\{1,2+\a_1+\a_2\}\l(\log\int_{A_{\fr{s}4,s}}|x|^{2\a_2}e^{u_2(x)}\mrm dV_g(x)-\fint_{B_s}u_2(x)\mrm dx-2(1+\a_2)\log s\r)\\
&\le&4\pi(1+\a_1)\l(\log\int_{B_{\fr{1}8}}|z|^{2\a_1}e^{v_1(z)}\mrm dz-\fint_{B_1}v_1(z)\mrm dx\r)\\
&+&2\pi\min\{1,2+\a_1+\a_2\}\l(\log\int_{A_{\fr{1}4,1}}e^{v_2(z)}\mrm dV_g(z)-\fint_{B_1}v_2(z)\mrm dz\r)+C\\
&\le&(1+\e)\int_{B_1}Q(v(z))\mrm dz+C\\
&=&(1+\e)\int_{B_s}Q(u(x))\mrm dx+C.
\eeqy
Finally, \eqref{inbdry} follows from Theorem \ref{mtdisk}:
\beqy
&&2\pi\sum_{i=1}^2\min\{1,2+\a_1+\a_2\}\l(\log\int_{A_{\fr{s}2,s}}|x|^{2\a_i}e^{u_i(x)}\mrm dx-\fint_{B_s}u_i(x)\mrm dx-2(1+\a_i)\log s\r)\\
&\le&2\pi\sum_{i=1}^2\min\{1,2+\a_1+\a_2\}\l(\log\int_{A_{\fr{1}2,1}}e^{v_i(z)}\mrm dz-\fint_{B_1}v_i(z)\mrm dz\r)+C\\
&\le&(1+\e)\int_{B_1}Q(v(z))\mrm dz+C\\
&=&(1+\e)\int_{B_s}Q(u(x))\mrm dx+C.
\eeqy
The final remark holds true because of the final remarks in Theorem \ref{mtdisk} and Lemma \ref{average}.
\epf

\blem\label{mtout}
For any $\e>0,\,\a_1,\a_2\in(-1,0],\,d>0$ small there exists $C=C_\e$ such that for any $p\in\Si,\,s\in\l(0,\fr{d}8\r)$ and $u\in H^1(\Si)^2$ with $u_i|_{\pa B_d(p)}\eq0$ one has
\beqa
\nonumber&&4\pi\sum_{i=1}^2(1+\a_{3-i})\log\int_{A_{2s,d}(p)}d(\cd,p)^{2\a_i}e^{u_i}\mrm dV_g+4\pi(1+\e)\sum_{i=1}^2(1+\a_i)\fint_{B_s(p)}u_i\mrm dV_g\\
\nonumber&+&8\pi(1+\e)\l((1+\a_1)^2+(1+\a_2)^2\r)\log s\\
\label{out}&\le&\int_{A_{s,d}(p)}Q(u)\mrm dV_g+\e\int_{B_d(p)}Q(u)\mrm dV_g+C,
\eeqa
\beqa
\nonumber&&4\pi(1+\a_2)\log\int_{A_{8s,d}(p)}d(\cd,p)^{2\a_1}e^{u_1}\mrm dV_g+4\pi(1+\e)(1+\a_1)\fint_{B_s(p)}u_1\mrm dV_g\\
\nonumber&+&2\pi\min\{1,2+\a_1+\a_2\}\l(\log\int_{A_{s,4s}(p)}d(\cd,p)^{2\a_2}e^{u_2}\mrm dV_g+(1+\e)\fint_{B_s(p)}u_2\mrm dV_g\r)\\
\nonumber&+&4\pi(1+\e)\l(2(1+\a_1)^2+\min\{1,2+\a_1+\a_2\}(1+\a_2)\r)\log s\\
\label{outmix}&\le&\int_{A_{s,d}(p)}Q(u)\mrm dV_g+\e\int_{B_d(p)}Q(u)\mrm dV_g+C,
\eeqa
\beqa
\nonumber&&2\pi\sum_{i=1}^2\min\{1,2+\a_1+\a_2\}\l(\log\int_{A_{s,2s}(p)}d(\cd,p)^{2\a_i}e^{u_i}\mrm dV_g+(1+\e)\fint_{B_s(p)}u_i\mrm dV_g\r)\\
\label{outbdry}&+&4\pi(1+\e)\min\l\{2+\a_1+\a_2,(2+\a_1+\a_2)^2\r\}\log s\\
\nonumber&\le&\int_{A_{s,d}(p)}Q(u)\mrm dV_g+\e\int_{B_d(p)}Q(u)\mrm dV_g+C.
\eeqa
The last statement holds true if $B_s(p)$ is replaced by a simply connected domain $\O_s$ belonging to $\mfrak A_{\d s}$ and such that $B_{\d s}(p)\sub\O_s\sub B_{\l(2+\fr{1}\d\r)s}(p)$ for some $\d>0$, with $C$ replaced by some $C_\d>0$.
\elem

\bpf
Just like Lemma \ref{mtin}, we will work with flat Euclidean balls, whose centers will be omitted. Moreover, it will not be restrictive to assume $\fint_{B_d}u_i(x)\mrm dx=0$ for both $i$'s.\\
Define, for $z\in B_d$ and $c_1,c_2\le-2(2+\a_1+\a_2)$,
$$v_i(z):=\l\{\bll(2c_i-c_{3-i})\log s&\tx{if }z\in B_s\\u_i\l(ds\fr{z}{|z|^2}\r)+(2c_i-c_{3-i})\log|z|&\tx{if }z\in A_{s,d}\earr\r.$$
By a change of variable we find, for $r\in\l\{\fr{1}8,\fr{1}4,\fr{1}2\r\}$,
\beqy
\int_{A_{s,rd}}|z|^{-4-2\a_i-2c_i+c_{3-i}}e^{v_i(z)}\mrm dz=\int_{B_{s,rd}}|z|^{-4-2\a_i}e^{u_i\l(ds\fr{z}{|z|^2}\r)}\mrm dz=(ds)^{-2-2\a_i}\int_{A_{\fr{s}r,d}}|x|^{2\a_i}e^{u_i(x)}\mrm dx\\
\int_{A_{rd,d}}e^{v_i(z)}\mrm dz\sim\int_{A_{rd,d}}|z|^{-4-2\a_i-2c_i+c_{3-i}}e^{v_i(z)}\mrm dz=(ds)^{-2-2\a_i}\int_{A_{s,\fr{s}r}}|x|^{2\a_i}e^{u_i(x)}\mrm dx.
\eeqy
Moreover, by Lemma \ref{average}, we get
\beqy
\l|\fint_{B_s}u_i(x)\mrm dx-\fint_{B_d}v_i(z)\mrm dz\r|&\le&\l|\fint_{B_s}u_i(x)\mrm dx-\fint_{\pa B_s}u_i(x)\mrm dx\r|+\l|\fint_{\pa B_s}u_i(x)\mrm dx-\fint_{\pa B_d}v_i(z)\mrm dz\r|\\
&+&\l|\fint_{B_d}v_i(z)\mrm dz-\fint_{B_d}v_i(z)\mrm dz\r|\\
&\le&C\sqrt{\int_{B_s}|\n u(x)|^2\mrm dx}+|(2c_i-c_{3-i})\log d|+C\sqrt{\int_{B_d}|\n v(z)|^2\mrm dz}\le\\
&\le&\e'\int_{B_s}Q(u(x))\mrm dx+\e'\int_{B_d}Q(v(z))\mrm dz+C_d.
\eeqy
Concerning the Dirichlet integral, we can write
\beqy
\int_{B_d}\n v_i(z)\cd\n v_j(z)\mrm dz&=&\int_{A_{s,d}}\l(\fr{d^2s^2}{|z|^4}\n u_i\l(sd\fr{z}{|z|^2}\r)\cd\n u_j\l(ds\fr{z}{|z|^2}\r)-(2c_i-c_{3-i})sd\fr{z}{|z|^2}\cd\n u_j\l(ds\fr{z}{|z|^2}\r)\r.\\
&-&\l.(2c_j-c_{3-j})sd\fr{z}{|z|^2}\cd\n u_i\l(sd\fr{z}{|z|^2}\r)+\fr{(2c_i-c_{3-i})(2c_j-c_{3-j})}{|z|^2}\r)\mrm dz\\
&=&\int_{A_{s,d}}\n u_i(x)\cd\n u_j(x)\mrm dx-(2c_i-c_{3-i})\int_{A_{s,d}}\fr{x}{|x|^2}\cd\n u_j(x)\mrm dx\\
&-&(2c_j-c_{3-j})\int_{A_{s,d}}\fr{x}{|x|^2}\cd\n u_i(x)\mrm dx-2\pi(2c_i-c_{3-i})(2c_j-c_{3-j})\log s\\
&+&2\pi(2c_i-c_{3-i})(2c_j-c_{3-j})\log d\\
&=&\int_{A_{s,d}}\n u_i(x)\cd\n u_j(x)\mrm dx+2\pi(2c_i-c_{3-i})\fint_{\pa B_s}u_j(x)\mrm dx\\
&+&2\pi(2c_j-c_{3-j})\fint_{\pa B_s}u_i(x)\mrm dx-2\pi(2c_i-c_{3-i})(2c_j-c_{3-j})\log s+C_d,
\eeqy
therefore, since $v$ has constant components in $B_s$,
$$\int_{B_d}Q(v(z))\mrm dz=\int_{A_{s,d}}Q(u(x))\mrm dx+2\pi\sum_{i=1}^2c_i\fint_{\pa B_s}u_i(x)\mrm dx-2\pi\l(c_1^2-c_1c_2+c_2^2\r)\log s+C.$$

The assertion of the lemma follows by applying Lemma \ref{mtloc} on $B_d$ to $v$ with different choices of $c_1,c_2$. If we take $c_1=c_2=-2(2+\a_1+\a_2)$, then we get
\beqy
&&4\pi\sum_{i=1}^2(1+\a_{3-i})\log\int_{A_{2s,d}}|x|^{2\a_i}e^{u_i(x)}\mrm dx\\
&\le&4\pi\sum_{i=1}^2(1+\a_{3-i})\l(\log\int_{B_{\fr{d}2}}|z|^{2\a_{3-i}}e^{v_i(z)}\mrm dz+2(1+\a_i)\log s\r)+C\\
&\le&(1+\e')\int_{B_d}Q(v(z))\mrm dz+4\pi\sum_{i=1}^2(1+\a_{3-i})\fint_{B_d}v_i(z)\mrm dz+16\pi(1+\a_1)(1+\a_2)\log s+C\\
&\le&(1+\e'')\int_{A_{s,d}}Q(u(x))\mrm dx+\e''\int_{B_d}Q(u(x))\mrm dx+4\pi\sum_{i=1}^2(1+\a_{3-i})\fint_{B_s}u_i(x)\mrm dx\\
&-&4\pi(1+\e'')(2+\a_1+\a_2)\sum_{i=2}^2\fint_{\pa B_s}u_i(x)\mrm dx+8\pi\l(2(1+\a_1)(1+\a_2)-(1+\e'')(2+\a_1+\a_2)^2\r)\log s+C\\
&\le&(1+\e''')\int_{B_s}Q(u(x))\mrm dx+\e'''\int_{B_d}Q(u(x))\mrm dx+4\pi\sum_{i=1}^2((1+\a_{3-i})-(1+\e''')(2+\a_1+\a_2))\fint_{B_s}u_i(x)\mrm dx\\
&+&8\pi\l(2(1+\a_1)(1+\a_2)-(1+\e''')(2+\a_1+\a_2)^2\r)\log s+C,
\eeqy
that is, re-naming $\e$ properly, \eqref{out}.\\

Choosing $c_1=-2(2+\a_1+\a_2)$ and $c_2=-2\min\{1,2+\a_1+\a_2\}=:-2m$, we get
\beqy
&&4\pi(1+\a_2)\log\int_{A_{8s,d}}|x|^{2\a_1}e^{u_1(x)}\mrm dx+2\pi m\log\int_{A_{s,4s}}|x|^{2\a_2}e^{u_2(x)}\mrm dx\\
&\le&4\pi(1+\a_2)\log\int_{B_{\fr{d}8}}|z|^{\max\{2+2\a_1+4\a_2,2\a_2\}}e^{v_1(z)}\mrm dz+2\pi m\log\int_{A_{\fr{d}4,d}}e^{v_2(z)}\mrm dz\\
&+&4\pi(2(1+\a_1)(1+\a_2)+m(1+\a_2))\log s+C\\
&\le&(1+\e')\int_{B_d}Q(v(z))\mrm dz+4\pi(1+\a_2)\fint_{B_d}v_1(z)\mrm dz+2\pi m\fint_{B_d}v_2(z)\mrm dz\\
&+&4\pi(1+\a_2)(2(1+\a_1)+m)\log s+C\\
&\le&(1+\e'')\int_{A_{s,d}}Q(u(x))\mrm dx+\e''\int_{B_d}Q(u(x))\mrm dx+4\pi(1+\a_2)\fint_{B_s}u_1(x)\mrm dx+2\pi m\fint_{B_s}u_2(x)\mrm dx\\
&-&4\pi(1+\e'')(2+\a_1+\a_2)\fint_{\pa B_s}u_1(x)\mrm dx-4\pi(1+\e'')m\fint_{\pa B_s}u_2(x)\mrm dx\\
&+&4\pi\l((1+\a_2)(2(1+\a_1)+m)-2(1+\e'')\l((2+\a_1+\a_2)^2-m(2+\a_1+\a_2)+m^2\r)\r)\log s+C\\
&\le&(1+\e''')\int_{A_{s,d}}Q(u(x))\mrm dx+\e'''\int_{B_d}Q(u(x))\mrm dx\\
&+&4\pi((1+\a_2)-(1+\e''')(2+\a_1+\a_2))\fint_{B_s}u_1(x)\mrm dx-2\pi(1+2\e''')m\fint_{B_s}u_2(x)\mrm dx\\
&+&4\pi((1+\a_2)(2(1+\a_1)+m)-2(1+\e''')((1+\a_1)(2+\a_1+\a_2)+(1+\a_2)m))\log s+C,
\eeqy
namely \eqref{outmix}.\\

Finally, taking $c_1=c_2=-2m$ one finds \eqref{outbdry}:
\beqy
&&2\pi\sum_{i=1}^2\min\{1,2+\a_1+\a_2\}\log\int_{A_{s,2s}}|x|^{2\a_i}e^{u_i(x)}\mrm dx\\
&\le&2\pi\sum_{i=1}^2m\l(\log\int_{A_{\fr{d}2,d}}e^{v_i(z)}\mrm dz+2(1+\a_i)\log s\r)+C\\
&\le&(1+\e')\int_{B_d}Q(v(z))\mrm dz+2\pi m\sum_{i=1}^2\fint_{B_d}v_i(z)\mrm dz+4\pi m(2+\a_1+\a_2)\log s+C\\
&\le&(1+\e'')\int_{B_d}Q(u(x))\mrm dx+\e''\int_{B_d}Q(u(x))\mrm dx+2\pi m\sum_{i=1}^2\fint_{B_s}u_i(x)\mrm dx\\
&-&4\pi(1+\e'')m\sum_{i=1}^2\fint_{\pa B_s}u_i(x)\mrm dx+4\pi\l(m(2+\a_1+\a_2)-2(1+\e'')m^2\r)\log s+C\\
&\le&(1+\e'')\int_{B_d}Q(u(x))\mrm dx+\e'''\int_{B_d}Q(u(x))\mrm dx-2\pi(1+2\e''')m\sum_{i=1}^2\fint_{B_s}u_i(x)\mrm dx\\
&-&4\pi(1+2\e''')m(2+\a_1+\a_2)\log s+C.
\eeqy
The final remark holds true, like in Lemma \ref{mtin}, because of the final remarks in Theorem \ref{mtdisk} and Lemma \ref{average}. In particular, when integrating by parts, one gets
$$\int_{B_\d\sm\O_s}\fr{x}{|x|^2}\cd\n u_i(x)\mrm dx=\int_{\pa\O_s}u_i(x)\ub{\fr{x}{|x|^2}\cd\nu(x)}_{=:f(x)}\mrm dx=,$$
with $\int_{\pa\O_s}f(x)\mrm dx=2\pi$ and $|f(x)|\le\fr{C}{s}\le\fr{C}{|\O_s|}$, therefore, by the Poincaré-Wirtinger inequality
\beqy
\l|\int_{\pa\O_s}u_i(x)f(x)\mrm dx-2\pi\int_{\pa\O_s}u_i(y)\mrm dy\r|&=&\l|\int_{\pa\O_s}f(x)\l(u_i(x)-\fint_{\pa\O_s}u_i(y)\mrm dy\r)\mrm dy\r|\\
&\le&C\fint_{\pa\O_s}\l|u_i(x)-\fint_{\pa\O_s}u(y)\mrm dy\r|\mrm dx\\
&\le&C\int_{\O_s}|\n u_i(x)|^2\mrm dx\\
&\le&\e\int_{\O_s}Q(u(x))\mrm dx+C_\e,
\eeqy
and $\l|\fint_{\pa\O_s}u_i(x)\mrm dx-\fint_{\pa B_s}u_i(x)\mrm dx\r|\le\e\int_{\O_s}Q(u(x))\mrm dx+C_\e$ by Lemma \ref{average}.
\epf

To prove Theorem \ref{impr} we also need the following lemma. It basically allows us to divide a disk in two domains in such a way that the integrals of two given functions are both split exactly in two.\\

\blem\label{split}
Consider $B:=B_1(0)\sub\R^2$ and $f_1,f_2\in L^1(B)$ such that $f_i>0$ a.e. $x\in B$ for both $i=1,2$ and $\int_Bf_1(x)\mrm dx=\int_Bf_2(x)\mrm dx=1$. Then, there exist $\t\in\S^1$ and $a\in(-1,1)$ such that
$$\int_{\{x\in B:\,x\cd\t<a\}}f_1(x)\mrm dx=\int_{\{x\in B:\,x\cd\t>a\}}f_2(x)\mrm dx=\fr{1}2$$
\elem

\bpf
Define, for $(\t,a)\in\S^1\x(-1,1)$,
$$I_1(\t,a):=\int_{\{x\in B:\,x\cd\t<a_1(\t)\}}f_1(x)\mrm dx.$$
For any given $\t$ there exists a unique $a_1(\t)$, smoothly depending on $\t$ such that $I_1(\t,a_1(\t))=\fr{1}2$. Define similarly $I_2(\t,a)$ and $a_2(\t)$.\\
Let us now show the existence of $\t$ such that $a_1(\t)=a_2(\t):=a$, hence the proof of the lemma will follow. Suppose by contradiction that $a_1(\t)<a_2(\t)$ for any $\theta$. Then, by definition, we get
$$a_1(-\t)=-a_1(\t)>-a_2(\t)=a_2(-\t),$$
which is a contradiction. One similarly excludes the case $a_1(\t)>a_2(\t)$.
\epf

\bpf[Proof of Theorem \ref{impr}]
From Lemma \ref{center} we have $\b\in\Si_{\rho_1,\ul\a_1}\cap\Si_{\rho_1,\ul\a_2},\si\in(0,\d)$ such that
$$\int_{B_\si(\b)}f_{1,u}\mrm dV_g\ge\d\q\q\q\int_{\Si\sm B_\si(\b)}f_{1,u}\mrm dV_g\ge\d\q\q\q\int_{B_\si(\b)}f_{2,u}\mrm dV_g\ge\d\q\q\q\int_{\Si\sm B_\si(\b)}f_{2,u}\mrm dV_g\ge\d.$$
Moreover, from Corollary \ref{scale}, we will suffice to prove the theorem for $\si\le2^{-\fr{6}\e-4}\d$.\\
We have to consider several cases, roughly following the proof of Proposition $3.2$ in \cite{mr13}.\\

\bite
\item[Case $1$]: $\int_{A_{\si,\d'}(\b)}f_{i,u}\mrm dV_g\ge\fr{\d}2$ for both $i=1,2$, where $\d':=2^{-\fr{3}\e}\d$.\\
As a first thing, we modify $u$ so that it vanishes outside $B_\d(\b)$: we take $n\in\l[1,\fr{2}\e\r]$ such that
$$\int_{A_{2^{n-1}\d',2^{n+1}\d'}(\b)}Q(u)\mrm dV_g\le\e\int_\Si Q(u)\mrm dV_g$$
and we define $u'_i$ as the solution of
$$\l\{\bll-\D u'_i=0&\tx{in }A_{2^{n-1}\d',2^{n+1}\d'}(\b)\\u'_i=u_i-\fint_{B_{2^n\d'}(\b)}u_i\mrm dV_g&\tx{on }\pa B_{2^n\d'}(\b)\\u'_i=0&\tx{on }\pa B_{2^{n+1}\d'}(\b)\earr\r.$$
$u'_i$ verifies, by Lemma \ref{trace},
$$\int_{A_{2^{n-1}\d',2^{n+1}\d'}(\b)}Q(u')\mrm dV_g\le C\int_{A_{2^{n-1}\d',2^{n+1}\d'}(\b)}Q(u)\mrm dV_g\le C\e\int_\Si Q(u)\mrm dV_g.$$
We obtained a function for which Lemma \ref{mtout} can be applied on $B_\d(\b)$. This was done at \emph{small price}, since the Dirichlet integral only increased by $\e$; moreover, $u'$ and $u$ coincide (up to an additive constant) on $B_{\d'}(\b)$, which is where both $f_{i,u}$'s attain most of their mass.\\
\bite
\item[Case $1.a$]: $\int_{B_\fr{\si}4(\b)}f_{i,u}\mrm dV_g\ge\fr{\d}2\tx{ for both }i=1,2.$\\
We apply Lemma \ref{mtin} to $u$ on $B_{\fr{\si}2(\b)}$, with $\a_i:=\a_{im}$ for $i=1,2$. From \eqref{in} we get
\beqa
\nonumber&&4\pi\sum_{i=1}^2(1+\a_{im})\l(\log\int_\Si\wt h_ie^{u_i}\mrm dV_g-\fint_{B_\fr{\si}2(\b)}u_i\mrm dV_g\r)-8\pi\l((1+\a_{1m})^2+(1+\a_{2m})^2\r)\log\fr{\si}2\\
\nonumber&\le&4\pi\sum_{i=1}^2(1+\a_{im})\l(\log\int_{B_{\fr{\si}4}(\b)}\wt h_ie^{u_i}\mrm dV_g-\fint_{B_\fr{\si}2(\b)}u_i\mrm dV_g\r)\\
\nonumber&-&8\pi\l((1+\a_{1m})^2+(1+\a_{2m})^2\r)\log\fr{\si}2+4\pi(2+\a_{1m}+\a_{2m})\log\fr{2}\d\\
\nonumber&\le&4\pi\sum_{i=1}^2(1+\a_{im})\l(\log\int_{B_{\fr{\si}4}(\b)}d(\cd,\b)^{2\a_{im}}e^{u_i}\mrm dV_g-\fint_{B_\fr{\si}2(\b)}u_i\mrm dV_g\r)\\
\nonumber&-&8\pi\l((1+\a_{1m})^2+(1+\a_{2m})^2\r)\log\fr{\si}2+C\\
\label{imprin}&\le&(1+\e)\int_{B_\fr{\si}2(\b)}Q(u)\mrm dV_g+C.
\eeqa
We then apply Lemma \ref{mtout} to $u'$ on $B_\d(\b)\sm B_\fr{\si}2(\b)$.
\beqa
\nonumber&&4\pi\sum_{i=1}^2(1+\a_{3-i,m})\log\int_\Si\wt h_ie^{u_i}\mrm dV_g+4\pi(1+\e)\sum_{i=1}^2(1+\a_{im})\fint_{B_\fr{\si}2(\b)}u_i\mrm dV_g\\
\nonumber&+&8\pi(1+\e)\l((1+\a_{1m})^2+(1+\a_{2m})^2\r)\log\fr{\si}2\\
\nonumber&\le&4\pi\sum_{i=1}^2(1+\a_{3-i,m})\log\int_{A_{\si,\d'}(\b)}\wt h_ie^{u'_i}\mrm dV_g+4\pi(1+\e)\sum_{i=1}^2(1+\a_{im})\fint_{B_\fr{\si}2(\b)}u'_i\mrm dV_g\\
\nonumber&+&8\pi(1+\e)\l((1+\a_{1m})^2+(1+\a_{2m})^2\r)\log\fr{\si}2+4\pi(2+\a_{1m}+\a_{2m})\log\fr{2}\d\\
\nonumber&\le&4\pi\sum_{i=1}^2(1+\a_{3-i,m})\log\int_{A_{\si,\d'}(\b)}d(\cd,\b)^{2\a_{im}}e^{u'_i}\mrm dV_g+4\pi(1+\e)\sum_{i=1}^2(1+\a_{im})\fint_{B_\fr{\si}2(\b)}u'_i\mrm dV_g\\
\nonumber&+&8\pi(1+\e)\l((1+\a_{1m})^2+(1+\a_{2m})^2\r)\log\fr{\si}2+C\\
\nonumber&\le&\int_{A_{\fr{\si}2,\d'}(\b)}Q(u')\mrm dV_g+\e\int_{B_{\d'}(\b)}Q(u')\mrm dV_g+C\\
\label{improut}&\le&\int_{A_{\fr{\si}2,\d'}(\b)}Q(u)\mrm dV_g+C\e\int_\Si Q(u)\mrm dV_g+C.
\eeqa
By summing \eqref{imprin} and \eqref{improut} and renaming properly $\e$ we get $J_{\rho_\e,\rho_\e}(u)\ge-L$ for $\rho_\e:=4\pi(2+\a_{1m}+\a_{2m})-\e$, which means, being $\e$ arbitrary, $J_\rho(u)\ge-L$.\\

\item[Case $1.b$]: $\int_{A_{4\si,\d'}(\b)}f_{i,u}\mrm dV_g\ge\fr{\d}4\tx{ for both }i=1,2.$\\
The result follows arguing as before, still applying Lemmas \ref{mtin}, \ref{mtout}, but this time on $B_{2\si}(\b)$ and $A_{2\si,\d'}(\b)$.\\

\item[Case $1.c$]:
$$\int_{B_\fr{\si}8(\b)}f_{1,u}\mrm dV_g\ge\fr{\d}2,\q\int_{A_{8\si,\d'}(\b)}f_{1,u}\mrm dV_g\ge\fr{\d}4,\q\int_{A_{\fr{\si}4,\si}(\b)}f_{2,u}\mrm dV_g\ge\fr{\d}2,\q\int_{A_{\si,4\si}(\b)}f_{2,u}\mrm dV_g\ge\fr{\d}4.$$
We still apply Lemmas \ref{mtin} and \ref{mtout}, respectively on $B_\si(\b)$ and $A_{\si,\d'}(\b)$, but this time we will exploit \eqref{inmix} and \eqref{outmix}: we get
\beqa
\nonumber&&4\pi(1+\a_{1m})\l(\log\int_\Si\wt h_1e^{u_1}\mrm dV_g-\fint_{B_\si(\b)}u_1\mrm dV_g\r)\\
\nonumber&+&2\pi\min\{1,2+\a_{1m}+\a_{2m}\}\l(\log\int_\Si\wt h_2e^{u_2}\mrm dV_g-\fint_{B_\si(\b)}u_2\mrm dV_g\r)\\
\nonumber&-&4\pi\l(2(1+\a_{1m})^2+\min\{1,2+\a_{1m}+\a_{2m}\}(1+\a_{2m})\r)\log\si+C\\
\nonumber&\le&4\pi(1+\a_{1m})\l(\log\int_{B_{\fr{\si}8}(\b)}d(\cd,\b)^{2\a_{1m}}e^{u_1}\mrm dV_g-\fint_{B_\si(\b)}u_1\mrm dV_g\r)\\
\nonumber&+&2\pi\min\{1,2+\a_{1m}+\a_{2m}\}\l(\log\int_{A_{\fr{\si}4,\si}(\b)}d(\cd,\b)^{2\a_{2m}}e^{u_2}\mrm dV_g-\fint_{B_\si(\b)}u_2\mrm dV_g\r)\\
\nonumber&-&4\pi\l(2(1+\a_{1m})^2+\min\{1,2+\a_{1m}+\a_{2m}\}(1+\a_{2m})\r)\log\si+C\\
\label{mixin}&\le&(1+\e)\int_{B_\si(\b)}Q(u)\mrm dV_g+C,
\eeqa
and
\beqa
\nonumber&&4\pi(1+\a_{2m})\log\int_\Si\wt h_1e^{u_1}\mrm dV_g+4\pi(1+\e)(1+\a_{1m})\fint_{B_\si(p)}u_1\mrm dV_g\\
\nonumber&+&2\pi\min\{1,2+\a_{1m}+\a_{2m}\}\l(\log\int_\Si\wt h_2e^{u_2}\mrm dV_g+(1+\e)\fint_{B_\si(p)}u_1\mrm dV_g\r)\\
\nonumber&+&4\pi(1+\e)\l(2(1+\a_{1m})^2+\min\{1,2+\a_{1m}+\a_{2m}\}(1+\a_{2m})\r)\log\si\\
\nonumber&\le&4\pi(1+\a_{2m})\log\int_{A_{8\si,d}(p)}d(\cd,p)^{2\a_{1m}}e^{u'_1}\mrm dV_g+4\pi(1+\e)(1+\a_{1m})\fint_{B_\si(p)}u'_1\mrm dV_g\\
\nonumber&+&2\pi\min\{1,2+\a_{1m}+\a_{2m}\}\l(\log\int_{A_{\si,4\si}(p)}d(\cd,p)^{2\a_{2m}}e^{u'_2}\mrm dV_g+(1+\e)\fint_{B_\si(p)}u'_2\mrm dV_g\r)\\
\nonumber&+&4\pi(1+\e)\l(2(1+\a_{1m})^2+\min\{1,2+\a_{1m}+\a_{2m}\}(1+\a_{2m})\r)\log\si+C\\
\nonumber&\le&\int_{A_{\si,\d'}(\b)}Q(u')\mrm dV_g+\e\int_{B_{\d'}(\b)}Q(u')\mrm dV_g+C,\\
\label{mixout}&\le&\int_{A_{\si,\d'}(\b)}Q(u)\mrm dV_g+C\e\int_\Si Q(u)\mrm dV_g+C,
\eeqa
As before, $J_\rho(u)\ge-L$ follows from \eqref{mixin}, \eqref{mixout} and a suitable redefinition of $\e$.\\

\item[Case $1.d$]:\\
$$\int_{A_{\fr{\si}4,\si}(\b)}f_{1,u}\mrm dV_g\ge\fr{\d}2,\q\int_{A_{\si,4\si}(\b)}f_{1,u}\mrm dV_g\ge\fr{\d}4,\q\int_{B_\fr{\si}8(\b)}f_{2,u}\mrm dV_g\ge\fr{\d}2,\q\int_{A_{8\si,\d'}(\b)}f_{2,u}\mrm dV_g\ge\fr{\d}4.$$
Here we argue as in case $1.b$, just exchanging the roles of $u_1$ and $u_2$.\\

\item[Case $1.e$]: $\int_{A_{\fr{\si}8,8\si}(\b)}f_{i,u}\mrm dV_g\ge\fr{\d}4$ for both $i=1,2$.\\
We would like to apply \eqref{inbdry} and \eqref{outbdry} and argue as in the previous cases. Anyway, we first need to define $\O_\si$ such that both components have some mass in both sets. We cover $A_{\fr{\si}8,8\si}(\b)$ with balls of radius $\fr{\si}{64}$; by compactness, we have $A_{\fr{\si}8,8\si}(\b)=\Cup_{l=1}^LB_\fr{\si}{64}(x_l)$, with $L$ not depending on $\si$, therefore there will be $x_{l_1},x_{l_2}$ such that
$\int_{B_\fr{\si}{64}\l(x_{l_i}\r)}f_{i,u}\mrm dV_g\ge\fr{\d}{4L}$.\\
We will proceed differently depending whether $x_{l_1}$ and $x_{l_2}$ are close or not.\\

\begin{figure}[h!]
\center
\includegraphics[width=0.4\linewidth]{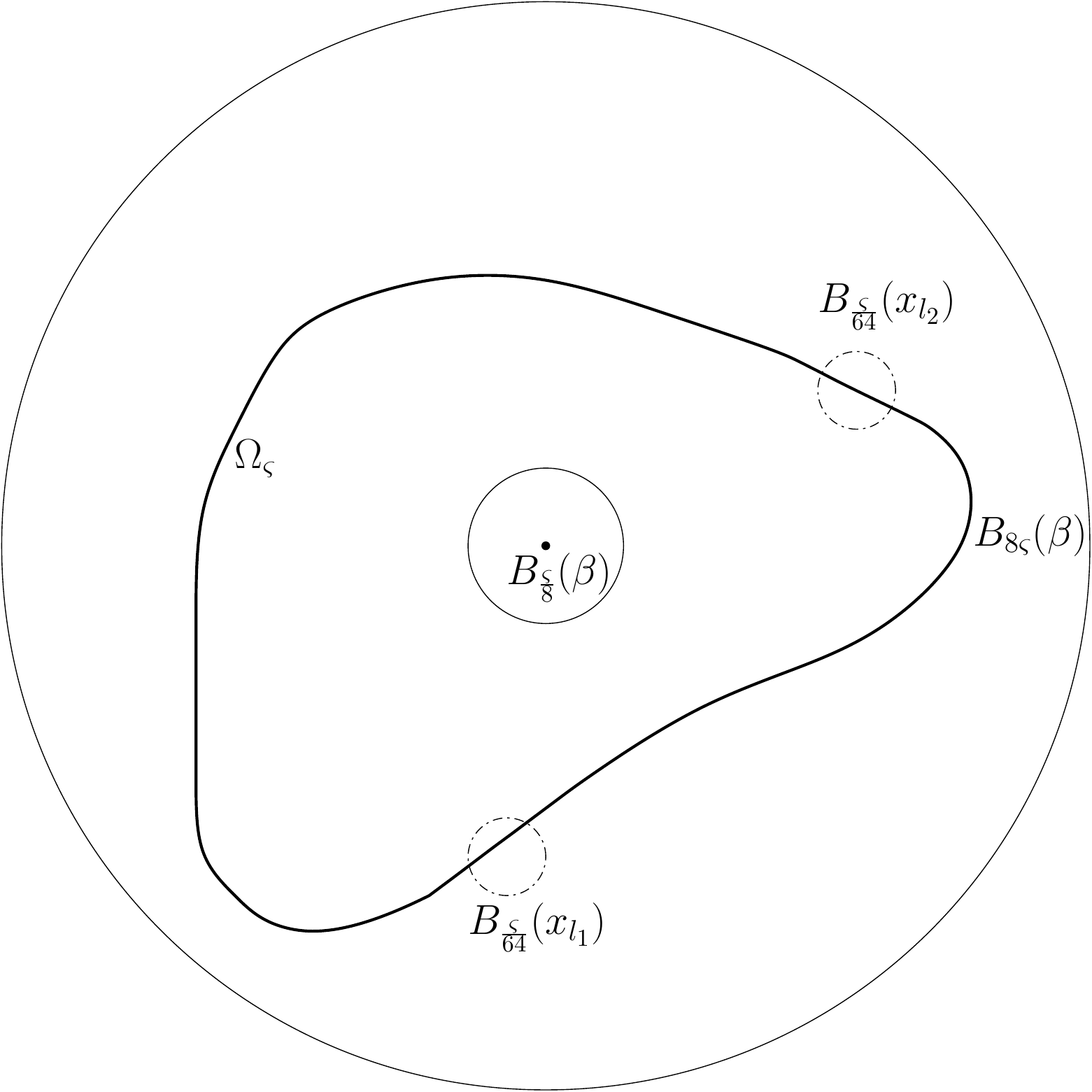}\q\q\q\q
\includegraphics[width=0.4\linewidth]{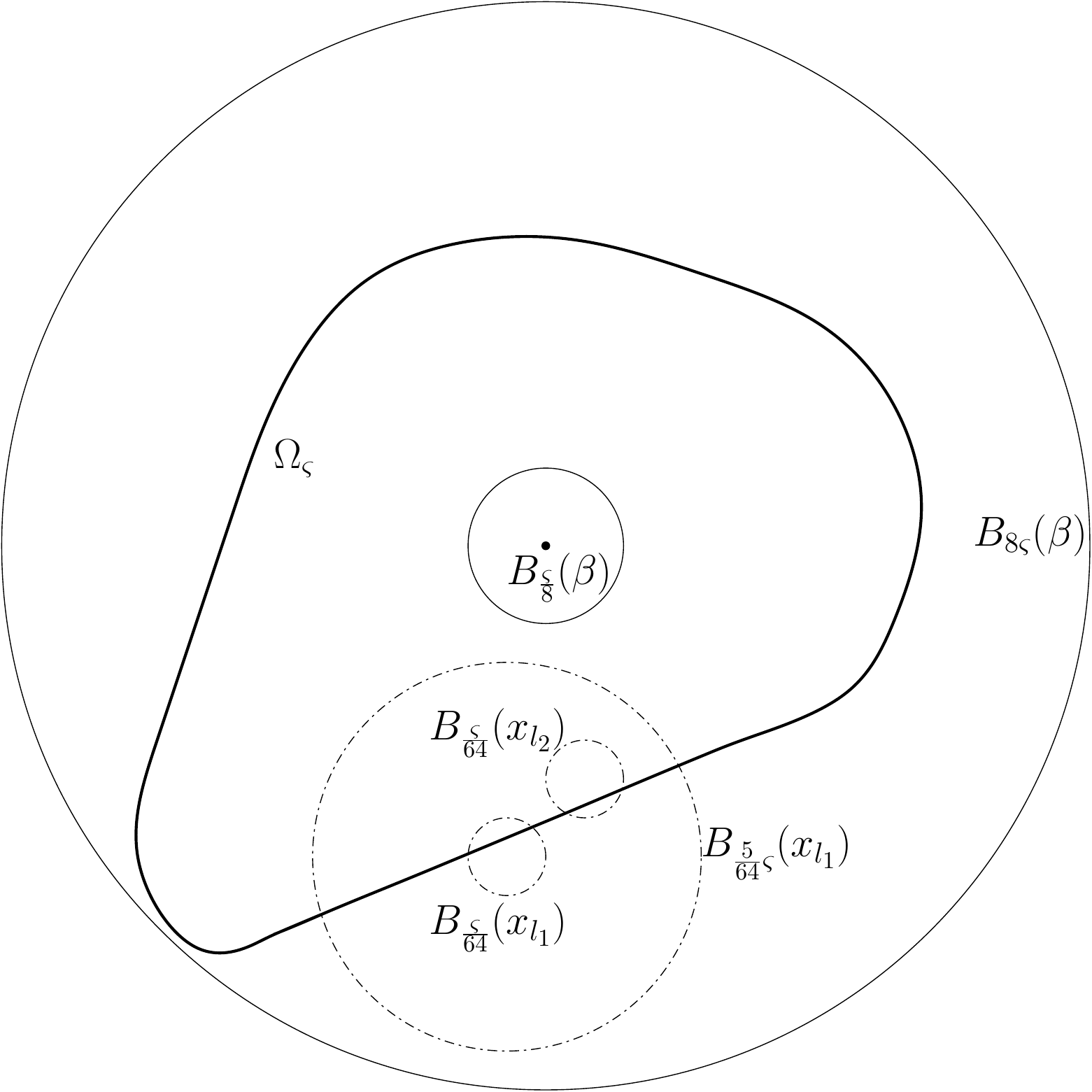}\\
\caption{The set $\O_\si$, respectively in the cases $1.e'$ and $1.e''$.}\label{figure}
\end{figure}

\bite
\item[Case $1.e'$]: $d(x_{l_1},x_{l_2})\ge\fr{\si}{16}$.\\
We divide each of the balls $B_\fr{\si}{64}(x_{l_1}),B_\fr{\si}{64}(x_{l_2})$ with a segment $\{x:(x-x_{l_i})\cd\t_i=a_i\}$, with $\t_i\in\S^1$ and $a_i\in\l(-\fr{\si}{64},\fr{\si}{64}\r)$, in such a way that
$$\int_{\l\{x\in B_\fr{\si}{64}\l(x_{l_i}\r),\l(x-x_{l_i}\r)\cd\t_i<a_i\r\}}f_{i,u}\mrm dV_g\ge\fr{\d}{8L}\q\int_{\l\{x\in B_\fr{\si}{64}\l(x_{l_i}\r),\l(x-x_{l_i}\r)\cd\t_i>a_i\r\}}f_{i,u}\mrm dV_g\ge\fr{\d}{8L}.$$
We can define $\O_\si$ as the region of $B_{\d'}(\b)$ delimited by the curve defined in the following way:\\
Since $d\l(B_\fr{\si}{32}(x_{l_1}),B_\fr{\si}{32}(x_{l_2})\r)\ge\fr{\si}{32}$, we can attach smoothly one endpoint of each segment without intersecting the two balls. We then join the other endpoint of each segment winding around $\b$.\\
Since $B_\fr{\si}{64}(x_{l_1})\sub A_{\fr{\si}{16},9\si}(\b)$, we can build $\O_\si$ in such a way that $\pa\O_\si\sub A_{\fr{\si}{32}(\b),10\si}$ and $\O_\si\in\mfrak A_{\d\si}$ (see \eqref{eq:ud} and Figure \ref{figure}). Moreover, by construction,
$$\int_{B_{\d'}(\b)\sm\O_\si}f_{i,u}\mrm dV_g\ge\fr{\d}{8L}\q\q\q\int_{\O_\si}f_{i,u}\mrm dV_g\ge\fr{\d}{8L},$$
hence Lemmas \ref{mtin} and \ref{mtout} still yield the proof.\\

\item[Case $1.e''$]: $d(x_{l_1},x_{l_2})\le\fr{\si}{16}$.\\
Since $B_\fr{\si}{64}(x_{l_1})\cup B_\fr{\si}{64}(x_{l_2})\sub B_{\fr{5}{64}\si}(x_{l_1})$, we apply Lemma \ref{split} to $f_i:=\fr{\wt h_ie^{u_i}}{\int_{B_{\fr{5}{64}\si}\l(x_{l_1}\r)}\wt h_ie^{u_i}\mrm dV_g}$ to find $\t\in\S^1,a\in\l(-\fr{5}{64}\si,\fr{5}{64}\si\r)$ such that
$$\int_{\l\{x\in B_{\fr{5}{64}\si}\l(x_{l_1}\r),\l(x-x_{l_1}\r)\cd\t<a\r\}}f_{i,u}\mrm dV_g\ge\fr{\d}{8L}\q\int_{\l\{x\in B_{\fr{5}{64}\si}\l(x_{l_1}\r),\l(x-x_{l_1}\r)\cd\t>a\r\}}f_{i,u}\mrm dV_g\ge\fr{\d}{8L}.$$
We now join smoothly (and without intersecting the balls) the endpoints of the segment $\{x:(x-x_{l_2})\cd\t=a\}$ with an arc winding around $\b$. Then, we define $\O_\si$ as the region of $B_\d'(\b)$ delimited by the curve made by such an arc and that segment.\\
Since $B_{\fr{5}{64}\si}(x_{l_1})\sub A_{\fr{3}{64}\si,9\si}(\b)$, as before we will have $B_\fr{\si}{32}(\b)\sub\O_\si\sub B_{10\si}(\b)$ and $\O_\si\in\mfrak A_{\d\si}$, and we can argue again as before because clearly
$$\int_{B_{\d'}(\b)\sm\O_\si}f_{i,u}\mrm dV_g\ge\fr{\d}{8L}\q\q\q\int_{\O_\si}f_{i,u}\mrm dV_g\ge\fr{\d}{8L}.$$

\eite
\eite
\item[Case $2$]: $\int_{\Si\sm B_{\d'}(\b)}f_{i,u}\mrm dV_g\ge\fr{\d}2$ for some $i$.\\
It will be not restrictive to assume $i=1$. If we also have $\int_{\Si\sm B_{\d''}(\b)}f_{2,u}\mrm dV_g\ge\fr{\d}2$, with $\d'':2^{-\fr{3}\e}\d'$, then we get $J_\rho(u)\ge-L$ by applying Lemma \ref{mtmacro}, as in the proof of Corollary \ref{scale}. Therefore we will assume
$$\int_{A_{\si,\d''}(\b)}f_{2,u}\mrm dV_g\ge\fr{\d}2.$$
The idea is to combine the previous arguments with a \emph{macroscopic} improved Moser-Trudinger inequality.\\
As a first thing, define $u''$ as the solution of
$$\l\{\bll-\D u''_i=0&\tx{in }A_{2^{n-1}\d'',2^{n+1}\d''}(\b)\\u''_i=u_i-\fint_{B_{2^n\d''}(\b)}u_i\mrm dV_g&\tx{on }\pa B_{2^n\d''}(\b)\\u''_i=0&\tx{on }\pa B_{2^{n+1}\d''}(\b)\earr\r.$$
with $n\in\l[1,\fr{2}\e\r]$ such that
$$\int_{A_{2^{n-1}\d'',2^{n+1}\d'}(\b)}Q(u'')\mrm dV_g\le C\int_{A_{2^{n-1}\d'',2^{n+1}\d''}(\b)}Q(u)\mrm dV_g\le C\e\int_\Si Q(u)\mrm dV_g.$$
Suppose $u$ satisfies the hypotheses of Case $1.a$, that is $\int_{A_{\si,\d''}(\b)}f_{i,u}\mrm dV_g\ge\fr{\d}2$ for both $i=1,2$. Then, clearly \eqref{imprin} still holds, whereas \eqref{improut} does not because we cannot estimate the integral of $\int_\Si\wt h_1e^{u_1}\mrm dV_g$ with the same integral evaluated over $A_{\si,\d''}$.\\
Anyway, by Jensen's inequality and Lemma \ref{average} we get
\beqy\log\int_{A_{\si,\d''}(\b)}\wt h_1e^{u'_1}\mrm dV_g&\ge&\log\int_{A_{\fr{\d''}2,\d''}(\b)}\wt h_1e^{u_1}\mrm dV_g-\fint_{B_{2^n\d''}(\b)}u_1\mrm dV_g\\
&\ge&\fint_{A_{\fr{\d''}2,\d''}(\b)}u_1\mrm dV_g+\log\l|A_{\fr{\d''}2,\d''}(\b)\r|+\fint_{A_{\fr{\d''}2,\d''}(\b)}\log\wt h_1\mrm dV_g-\fint_{B_{2^n\d''}(\b)}u_1\mrm dV_g\\
&\ge&-\e\int_\Si Q(u)\mrm dV_g-C,
\eeqy
hence we obtain
\beqa
\nonumber&&4\pi(1+\a_{2m})\fint_{B_{2^n\d''}(\b)}u_1\mrm dV_g+4\pi(1+\a_{1m})\log\int_\Si\wt h_2e^{u_2}\mrm dV_g\\
\nonumber&+&4\pi(1+\e)\sum_{i=1}^2(1+\a_{im})\fint_{B_\fr{\si}2(\b)}u_i\mrm dV_g+8\pi(1+\e)\l((1+\a_{1m})^2+(1+\a_{2m})^2\r)\log\fr{\si}2\\
\label{improut2}&\le&\int_{A_{\fr{\si}2,\d''}(\b)}Q(u)\mrm dV_g+C\e\int_\Si Q(u)\mrm dV_g+C.
\eeqa
Now, by Jensen's inequality and a variation of the \emph{localized} Moser-Trudinger inequality \eqref{loc},
\beqa
\nonumber&&4\pi\l(1+\min_{m'\ne m}\a_{1m'}\r)\l(\log\int_{\Si\sm B_{\d'}(\b)}\wt h_1e^{u_1}\mrm dV_g-\int_\Si u_1\mrm dV_g\r)\\
\nonumber&\le&4\pi\sum_{i=1}^2\l(1,1+\min_{m'\ne m}\a_{im'}\r)\l(\log\int_{\Si\sm B_{\d'}(\b)}\wt h_ie^{u_i}\mrm dV_g-\int_\Si u_i\mrm dV_g\r)+C\\
\label{locout}&\le&(1+\e)\int_{\Si\sm B_\fr{\d'}2(\b)}Q(u)\mrm dV_g+C.
\eeqa
By summing \eqref{imprin}, \eqref{improut2} and \eqref{locout} we get $J_{\rho_{1\e},\rho_{2\e}}(u)\ge-L$, with
$$\rho_{1\e}:=4\pi\min\l\{2+\a_{1m}+\a_{2m},1+\a_{1m}+\min_{m'\ne m}(1+\a_{1m'})\r\}-\e\q\rho_{2\e}:=4\pi(2+\a_{1m}+\a_{2m})-\e,$$
therefore $J_\rho(u)\ge-L$. We argue similarly if we are under the condition of Cases $1.b,\,1.c,\,1.d,\,1.e$.
\eite
The proof is thereby concluded.
\epf

\section{Proof of Theorem \ref{ex}}\label{s:ex}

We are finally in position to prove the main existence theorem of this paper. Its proof will follow by showing that low sub-levels are \emph{dominated} by the space $\mcal X$ (see \cite{hat}, page 528), which is not contractible by the results contained in Section \ref{s:top}. In particular, we have the following lemma, whose proof is given below.\\

\blem\label{homo}
For $L\gg0$ large enough there exist maps $\Phi:\mcal X\to J_\rho^{-L}$ and $\Psi:J_\rho^{-L}\to\mcal X$ such that $\Psi\c\Phi$ is homotopically equivalent to $\mrm{Id}_{\mcal X}$.
\elem

\bpf[Proof of Theorem \ref{ex}]
Suppose by contradiction that the system \eqref{toda} has no solutions. By Lemma \ref{deform}, $J_\rho^{-L}$ is a deformation retract of $J_\rho^L$, hence by Corollary \ref{contr} it is contractible. Let $H(\z,s):\mcal X\x[0,1]\to\mcal X$ be the homotopy equivalence defined in Lemma \ref{homo} and let $H'$ be a homotopy equivalence between a constant map and $\mrm{Id}_{J_\rho^{-L}}$.\\
Then $H''(\z,s)=\Psi(H'(\Phi(\z),s)):\mcal X\x[0,1]\to\mcal X$ is an equivalence between the maps $\Psi\c\Phi$ and a constant and $H''\ast H$ is an equivalence between $\mrm{Id}_{\mcal X}$ and a constant map. This means that $\mcal X$ is contractible, in contradiction with Theorem \ref{top}.
\epf

To prove Lemma \ref{homo} we need the following estimate. Notice that the choice of $\tau$ (see the proof of Lemma \ref{center}), which was not relevant in all the rest of this paper, will be made in the proof of this lemma to let the following result hold true.\\

\blem\label{estim}
Let $\d$ be as in Lemma \ref{center}, $\b_i(u),\si_i(u)$ be as in Corollary \ref{scale} and $\Phi^\la$ as in Theorem \ref{test}. Then, for a suitable choice of $\tau$, there exists $C_0>0,\d'\in(0,\d)$ such that:
\bite
\item If either $t\ge1-\fr{C_0}\la$ or $\l\{\bl t>\fr{1}2\\x_1=x_2=:p_m\\\rho_1,\rho_2<4\pi(2+\a_{1m}+\a_{2m})\earr\r.$, \q\q\q then $\si_1\l(\Phi^\la(\z)\r)\ge\d'$;\\
otherwise, $\si_1\l(\Phi^\la(\z)\r)<\d$ and $\b_1\l(\Phi^\la(\z)\r)=x_1$.
\item If either $t\le\fr{C_0}\la$ or $\l\{\bl t<\fr{1}2\\x_1=x_2=:p_m\\\rho_1,\rho_2<4\pi(2+\a_{1m}+\a_{2m})\earr\r.$, \q\q\q then $\si_2\l(\Phi^\la(\z)\r)\ge\d'$;\\
otherwise, $\si_2\l(\Phi^\la(\z)\r)<\d$ and $\b_2\l(\Phi^\la(\z)\r)=x_2$.
\eite
\elem

\bpf
We will only prove the statements involving $\si_1$ and $f_{1,\Phi^\la(\z)}$, since the same proof will work for the rest, up to switching indexes $i=1,2$. We will show the proof only in the case $x_2=p_m',\,\rho_2>4\pi(2+\a_{1m'}+\a_{2m'})$, which is somehow trickier because $\ph_1$ does not vanish when $t\ge1-\fr{1}\la$. Let us write
\beqy\ph_1&=&\l(\ph_{1,p_m}^{\la(1-t)}+\ph_{1,p_{m'}}^{\la t}\r)=\l(\ph_{1,p_m}^{\la(1-t)}-2\log\max\l\{1,(\la t)^{2(2+\a_{1m'}+\a_{2m'})}d(\cd,p_{m'})^{2(1+\a_{1m'})}\r\}\r),\\
\ph_2&=&\l(\ph_{2,p_m}^{\la(1-t)}+\ph_{2,p_{m'}}^{\la t}\r)=\l(\ph_{2,p_m}^{\la(1-t)}-2\log\max\l\{1,(\la td(\cd,p_{m'}))^{2(2+\a_{1m'}+\a_{2m'})}\r\}\r).
\eeqy
From the definition of $\si_1$, we have to show that, if $t\ge1-\fr{C_0}\la$, then
$$\int_{B_{\d'}(p_{m''})}f_{1,\Phi^\la(\z)}\mrm dV_g<\tau\q\q\q\fa\,m''=1,\ds,M.$$
It is not hard to see that, for any $m''\ne m,m'$,
$$\int_{B_{\d'}(p_{m''})}f_{1,\Phi^\la(\z)}\mrm dV_g\le C'{\d'}^{2(1+\a_{1m''})},$$
which is smaller than any given $\tau$ if $\d'$ is taken small enough. Roughly speaking, $f_{1,\Phi^\la(\z)}$ cannot attain mass too near $p_m$ because its scale depends on $\la(1-t)$ which is bounded from above. Moreover, $\ph_{1,p_m}^{\la(1-t)}$ is constant in $B_{(\la(1-t))^{-\fr{2+\a_{1m}+\a_{2m}}{1+\a_{1m}}}(p_m)}$, hence for large $C_0$
$$\int_{B_{C_0^{-1-\fr{(2+\a_{1m}+\a_{2m})}{1+\a_{1m}}}}(p_m)}f_{1,\Phi^\la(\z)}\mrm dV_g\le CC_0^{2(1+\a_{1m})}\int_{B_{C_0^{-\fr{(2+\a_{1m}+\a_{2m})}{1+\a_{1m}}}}(p_m)}d(\cd,p_m)^{2\a_{1m}}\mrm dV_g\le\fr{1}2<\tau.$$
On the other hand, a part of the mass of $f_{1,\Phi^\la(\z)}$ could actually concentrate near $p'_m$, but not all of it. Here, we will have to take $\tau$ properly. Since
\beqy\int_{B_\fr{\d}2(p_{m'})}\wt h_1e^{\ph_1-\fr{\ph_2}2}\mrm dV_g&\le&Ce^{\int_\Si\l(\ph_{1,p_m}^{\la(1-t)}-\fr{\ph_{2,p_m}^{\la(1-t)}}2\r)\mrm dV_g}\l(\int_{B_{(\la t)^{-\fr{2+\a_{1m'}+\a_{2m'}}{1+\a_{1m'}}}}(p_m')}d(\cd,p_{m'})^{2\a_{1m'}}\mrm dV_g\r.\\
&+&(\la t)^{-4\l(2+\a_{1m'}+\a_{2m'}\r)}\int_{A_{(\la t)^{-\fr{2+\a_{1m'}+\a_{2m'}}{1+\a_{1m'}}},\fr{1}{\la t}}(p_{m'})}d(\cd,p_m)^{-2(2+\a_{1m'})}\mrm dV_g\\
&+&\l.(\la t)^{-4\l(2+\a_{1m'}+\a_{2m'}\r)}\int_{A_{\fr{1}{\la t},\fr{\d}2}(p_{m'})}d(\cd,p_{m'})^{2\a_{2m'}}\mrm dV_g\r)\\
&\le&Ce^{\int_\Si\l(\ph_{1,p_m}^{\la(1-t)}-\fr{\ph_{2,p_m}^{\la(1-t)}}2\r)\mrm dV_g}(\la t)^{-4\l(2+\a_{1m'}+\a_{2m'}\r)},
\eeqy
and
\beqy
\int_{A_{\fr{\d}2,\d}(p_{m'})}\wt h_1e^{\ph_1-\fr{\ph_2}2}\mrm dV_g&\ge&\fr{1}Ce^{\int_\Si\l(\ph_{1,p_m}^{\la(1-t)}-\fr{\ph_{2,p_m}^{\la(1-t)}}2\r)\mrm dV_g}(\la t)^{-4\l(2+\a_{1m'}+\a_{2m'}\r)}\int_{A_{\fr{\d}2,\d}(p_{m'})}d(\cd,p_{m'})^{2\a_{2m'}}\mrm dV_g\\
&\ge&\fr{1}Ce^{\int_\Si\l(\ph_{1,p_m}^{\la(1-t)}-\fr{\ph_{2,p_m}^{\la(1-t)}}2\r)\mrm dV_g}(\la t)^{-4\l(2+\a_{1m'}+\a_{2m'}\r)},
\eeqy
then
$$\int_{B_\fr{\d}2(p_{m'})}f_{1,\Phi^\la(\z)}\mrm dV_g<\fr{\int_{B_\fr{\d}2(p_{m'})}f_{1,\Phi^\la(\z)}\mrm dV_g}{\int_{B_\d(p_{m'})}f_{1,\Phi^\la(\z)}\mrm dV_g}\le\fr{C^2}{1+C^2}.$$
Therefore, setting $\tau:=\fr{C^2}{1+C^2}$, we proved the first part of the Lemma.\\

Let us now assume $t\le1-\fr{C_0}\la$. From the proof of Lemma \ref{logint}, we deduce that the ratio $\fr{\int_{B_\d(p_m)}f_{1,\Phi^\la(\z)}\mrm dV_g}{\int_{B_\d(p_{m''})}f_{1,\Phi^\la(\z)}\mrm dV_g}$ increases arbitrarily as $\la(1-t)$ increases. Therefore, for large $C_0$, most of the mass of $f_{1,\Phi^\la(\z)}$ will be around $p_m$, hence by definition we will have $\b_1\l(\Phi^\la(\z)\r)=p_m$ and $\si_1\l(\Phi^\la(\z)\r)<\d$.
\epf

\bpf[Proof of Lemma \ref{homo}]
Let $\d$ be as in Lemma \ref{center}, $\b_i(u),\si_i(u)$ be as in Corollary \ref{scale} and $\d'$ be as in Lemma \ref{estim}. Take now $L$ so large that Corollary \ref{scale} and Theorem \ref{impr} apply.\\
We define $\Phi=\Phi^{\la_0}$ as in Theorem \ref{test}, with $\la_0$ such that $\Phi^\la(\mcal X)\sub J_\rho^{-2L}$ for any $\la\ge\la_0$. As for $\Psi:J_\rho^{-2L}\to\mcal X$, we write
$$\Psi(u)=\l(\b_1(u),\b_2(u),t'(\si_1(u),\si_2(u))\r)\q\q\q\tx{with }t'(\si_1(u),\si_2(u))=\l\{\bll0&\tx{if }\si_2(u)\ge\d'\\\fr{\d'-\si_2(u)}{2\d'-\si_1(u)-\si_2(u)}&\tx{if }\si_1(u),\si_2(u)\le\d'\\1&\tx{if }\si_1(u)\ge\d'\earr\r..$$
Let us verify the well-posedness of $\Psi$. The definition of $t'$ makes sense because, from Corollary \ref{scale}, $J_\rho(u)<-L$ implies $\min\{\si_1(u),\si_2(u)\}\le\d'$. Moreover, if $t'>0$ (respectively, $t'<1$), then $\si_1<\d$ is well-defined (respectively, $\si_2<\d$ is well-defined), hence $\b_1$ (respectively, $\b_2$) is also defined. Finally, $\Psi$ is mapped on $\mcal X$ because, from Theorem \ref{impr}, when $J_\rho(u)<-L$ we cannot have $(\b_1(u),\b_2(u),t'(\si_1(u),\si_2(u)))=\l(p_m,p_m,\fr{1}2\r)$ with $\rho_1,\rho_2<4\pi(2+\a_{1m}+\a_{2m})$.\\

To get a homotopy between the two maps, we first let $\la$ tend to $+\infty$, in order to get $x_1$ and $x_2$, then we apply a linear interpolation for the parameter $t$.\\
Writing $\Psi\l(\Phi^\la(\z)\r)=\l(\b_1^\la(\z),\b_2^\la(\z),t'^\la(\z)\r)$, we have $F=F_2\ast F_1$, with
$$F_1:(\z,s)=((x_1,x_2,t),s)\mapsto\l(\b_1^\fr{\la_0}{1-s}(\z),\b_2^\fr{\la_0}{1-s}(\z),t'^{\la_0}(\z)\r)$$
$$F_2:\l(x_1,x_2,t'^{\la_0}(\z)\r)\mapsto\l(x_1,x_2,(1-s)t'^{\la_0}(\z)+st\r).$$
We have to verify that all is well-defined.\\
If we cannot define $\b_1^{\fr{\la_0}{1-s}}(\z)$, then by Lemma \ref{estim} we either have $t\ge1-\fr{C_0(1-s)}{\la_0}\ge1-\fr{C_0}{\la_0}$ or we are on the first half of the punctured segment. By the same lemma, we get $\si_1\l(\Phi^{\la_0}(\z)\r)\ge\d'$,that is $t'^{\la_0}(\z)=1$. For the same reason, if $\b_2^\fr{\la_0}{1-s}(\z)$ is not defined, then $t'^{\la_0}(\z)=0$, so $F_1:\mcal X\x[0,1]\to\Si_{\rho_1,\ul\a_1}\star\Si_{\rho_2,\ul\a_2}$ makes sense.\\
Its image is actually contained in $\mcal X$ because, from Lemma \ref{estim}, if $x_1=x_2$ and $\rho<4\pi\l(\o_{\ul\a_1}(x)+\o_{\ul\a_2}(x)\r)$, then either $t'^{\la_0}(\z)\in\{0,1\}$, hence in particular it does not equal $\fr{1}2$.\\
Concerning $F_2$, the previous lemma implies $\b_1^\fr{\la_0}{1-s}(\z)=x_1$ if $t\le1-\fr{C_0}\la(1-s)$, hence in particular passing to the limit as $s\to1$, if $t<1$. A similar condition holds for $\b_2$, which gives $F_2(\cd,0)=F_1(\cd,1)$. If $x_1$ is not defined then $t'^{\la_0}(\z)=1$, hence $(1-s)t'^{\la_0}(\z)+st=1$, and similarly there are no issues when $x_2$ cannot be defined. Finally, by the argument used before, if $x_1=x_2=p_m$ and $\rho_1,\rho_2<4\pi(2+\a_{1m}+\a_{2m})$, then $(1-s)t'^{\la_0}(\z)+st\ne\fr{1}2$.
\epf

\section{Proofs on the non-existence results}\label{s:nonex}

\subsection{Proof of Theorems \ref{nonex1},\ref{nonex2}}

In this section we will consider some cases that are not covered by Theorem \ref{ex}. They are both inspired by \cite{barmal} (Propositions $5.7$ and $5.8$, respectively).\\

We start by considering the case of the unit disk $\l(\B,g_0\r)$ with one singularity in its center. Even though we are not dealing with a closed surface, most of the variational theory for the Liouville equations and the Toda system can be applied in the very same way to Euclidean domains (or surfaces with boundary) with Dirichlet boundary conditions. This was explicitly pointed out in \cite{barmal,bat2} for the Liouville equations, but still holds true for the Toda system, since blow-up on the boundary was excluded in \cite{lwzhao}. In particular, the general existence result contained in \cite{bjmr} holds on any non-simply connected open domain of the plane, since such domains can be retracted on a \emph{bouquet} of circles.\\
Concerning simply connected domains, we have minimizing solutions in the range of parameters $\rho_1<4\pi(1+\a_1),\rho_2<4\pi(1+\a_2)$, as well as the configuration $(M_1,M_2,M_3)=(1,1,0)$ in Theorem \ref{ex}. The region generating minimizing solution is colored in orange in Figure $4$, the region generating min-max solutions in colored in green.\\
By Poho\v zaev identity we show that most of the remaining set of parameters yields no solutions, colored in blue in the figure, and this holds in particular if one or both $\rho_i$'s are large enough.\\

Theorem \ref{nonex1} still holds if $\a_1=\a_2=0$, that is if we consider the regular Toda system. Here we still have solutions in the second square $(4\pi,8\pi)^2$; arguing as in \cite{mr13} we get low sub-levels being dominated by a space which is homeomorphic to $\R^6\sm\R^3\simeq\S^2$. This was confirmed in \cite{lwyang}, where the degree for the Toda system is computed, and in this case it equals $-1$. Figure $4$ shows that there might not be solutions in each of all the other squared which are delimited by integer numbers of $4\pi$. In particular, this shows that the degree is $0$ in all these regions.\\

\begin{figure}[h!]
\center
\includegraphics[width=0.3\linewidth]{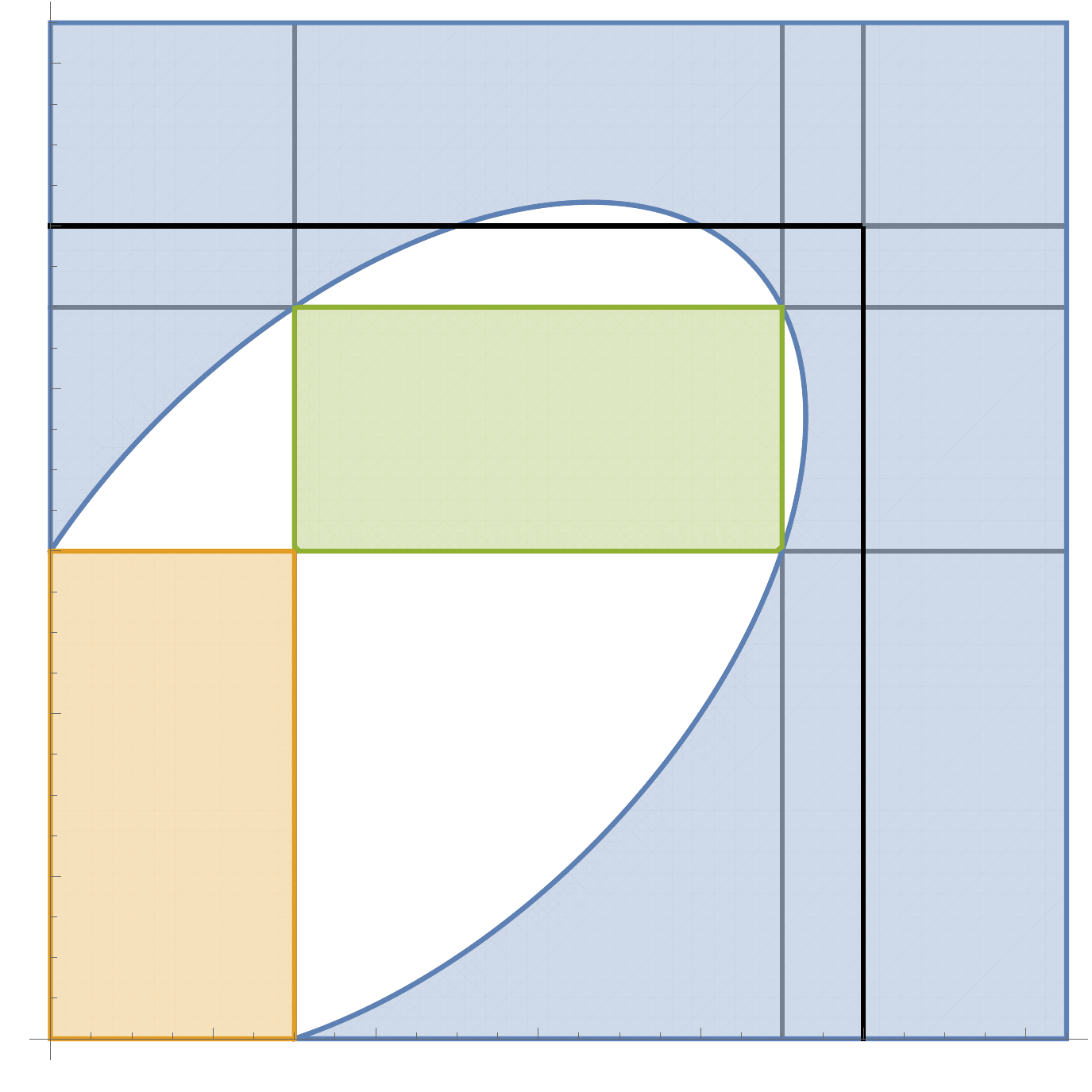}\\
\caption{Values of $\rho$ which yield existence and non-existence results for $\l(\B^2,g_0\r)$.}
\end{figure}

\bpf[Proof of Theorem \ref{nonex1}]
Let $u=(u_1,u_2)$ be a solution of \eqref{todadisk}. Since both components vanish on the boundary, for any $x\in\pa\B$ one has $\n u_i(x)=(\n u_i(x)\cd\nu(x))\nu(x)=:\pa_\nu u_i(x)\nu(x)$ for both $i=1,2$.\\
Therefore, one can apply a standard Poho\v zaev identity:
\beqy
&&\int_{\pa\B}\l((\pa_\nu u_1)^2+\pa_\nu u_1\pa_\nu u_2+(\pa_\nu u_2)^2\r)\mrm d\s\\
&=&2\int_{\pa\B}\l((\pa_\nu u_1)^2-\fr{|\n u_1|^2}2+\pa_\nu u_1\pa_\nu u_2-\fr{\n u_1\cd\n u_2}2+(\pa_\nu u_2)^2-\fr{|\n u_2|^2}2\r)\mrm d\s\\
&=&\int_\B\l(2(x\cd\n u_1(x))\D u_1(x)+(x\cd\n u_1(x))\D u_2(x)+(x\cd\n u_2(x))\D u_1(x)+2(x\cd\n u_2(x))\D u_2(x)\r)\mrm dx\\
&=&-3\fr{\rho_1}{\int_\B|x|^{2\a_1}e^{u_1(x)}\mrm dx}\int_\B(x\cd\n u_1(x))|x|^{2\a_1}e^{u_1(x)}\mrm dx-3\fr{\rho_2}{\int_\B|x|^{2\a_2}e^{u_2(x)}\mrm dx}\int_\B(x\cd\n u_2(x))|x|^{2\a_2}e^{u_2(x)}\mrm dx\\
&=&6\rho_1\fr{\int_{\pa\B}|\cd|^{2\a_1}e^{u_1}\mrm d\s}{\int_\B|x|^{2\a_1}e^{u_1(x)}\mrm dx}+6(1+\a_1)\rho_1+6\rho_2\fr{\int_{\pa\B}|\cd|^{2\a_1}e^{u_2}\mrm d\s}{\int_\B|x|^{2\a_2}e^{u_2(x)}\mrm dx}+6(1+\a_2)\rho_2.
\eeqy
For the boundary integral, we can perform an algebraic manipulation, use H\"older's inequality and then integrate by parts:
\beqy
&&\int_{\pa\B}\l((\pa_\nu u_1)^2+\pa_\nu u_1\pa_\nu u_2+(\pa_\nu u_2)^2\r)\mrm d\s\\
&=&\int_{\pa\B}\l(\fr{1}4(\pa_\nu u_1+2\pa_\nu u_2)^2+\fr{3}4(\pa_\nu u_1)^2\r)\mrm d\s\\
&\ge&\fr{1}{2\pi}\l(\fr{1}4\l(\int_{\pa\B}\pa_\nu u_1\mrm d\s+2\int_{\pa\B}\pa_\nu u_2\mrm d\s\r)^2+\fr{3}4\l(\int_{\pa\B}\pa_\nu u_1\mrm d\s\r)^2\r)\\
&=&\fr{1}{2\pi}\l(\l(\int_{\pa\B}\pa_\nu u_1\mrm d\s\r)^2+\l(\int_{\pa\B}\pa_\nu u_1\mrm d\s\r)\l(\int_{\pa\B}\pa_\nu u_2\mrm d\s\r)+\l(\int_{\pa\B}\pa_\nu u_2\mrm d\s\r)^2\r)\\
&=&\fr{1}{2\pi}\l(\l(\int_\B\D u_1(x)\mrm dx\r)^2+\l(\int_\B\D u_1(x)\mrm dx\r)\l(\int_\B\D u_2(x)\mrm dx\r)+\l(\int_\B\D u_2(x)\mrm dx\r)^2\r)\\
&=&\fr{3}{2\pi}\l(\rho_1^2-\rho_1\rho_2+\rho_2^2\r).
\eeqy
Therefore, we get as a necessary condition for existence of solutions:
\beqy
\rho_1^2-\rho_1\rho_2+\rho_2^2&=&4\pi\l(\rho_1\fr{\int_{\pa\B}|\cd|^{2\a_1}e^{u_1}\mrm d\s}{\int_\B|x|^{2\a_1}e^{u_1(x)}\mrm dx}+(1+\a_1)\rho_1+\rho_2\fr{\int_{\pa\B}|\cd|^{2\a_1}e^{u_2}\mrm d\s}{\int_\B|x|^{2\a_2}e^{u_2(x)}\mrm dx}+(1+\a_2)\rho_2\r)\\
&>&4\pi(1+\a_1)\rho_1+4\pi(1+\a_2)\rho_2.
\eeqy
This concludes the proof.
\epf

Let us now consider the standard sphere $\l(\S^2,g_0\r)$ with two antipodal singularities.\\

In Theorem \ref{nonex2} we perform a stereographic projection that transforms the solutions of \eqref{toda} on $\S^2$ on entire solutions on the plane, and then we use a Poho\v zaev identity for the latter problem, getting necessary algebraic condition for the existence of solutions.\\
Such a Poho\v zaev identity yields an algebraic condition which is similar to the one which appears in Theorem \ref{nonex1}. It can be deduced in the same way as was done in \cite{cl91} for the scalar Liouville equation.\\

\bthm\label{poho2}
Let $H_1,H_2\in L^p_{\mrm{loc}}\l(\R^2\r)$ be such that, for suitable $a>0,\,b>-2,\,C>0$,
$$\fr{|x|^a}C\le H_i(x)\le C|x|^b\q\fa\,x\in B_1(0)\sm\{0\}\q\q\q0<H_i(x)\le C|x|^a\q\fa\,x\in\R^2\sm B_1(0);$$
let $U=(U_1,U_2)$ be a solution of
$$\l\{\bll-\D U_1=2H_1e^{U_1}-H_2e^{U_2}&\tx{on }\R^2\\-\D U_2=2H_2e^{U_2}-H_1e^{U_1}&\tx{on }\R^2\\\int_{\R^2}\l(|x|^a+|x|^b\r)e^{U_1(x)}\mrm dx<+\infty\\\int_{\R^2}\l(|x|^a+|x|^b\r)e^{U_2(x)}\mrm dx<+\infty\earr\r.$$
and define
$$\rho_i:=\int_{\R^2}H_i(x)e^{U_i(x)}\mrm dx,\q\q\q,\tau_i:=\int_{\R^2}(x\cd\n H_i(x))e^{U_i(x)}\mrm dx,\q\q\q i=1,2.$$
Then,
$$\rho_1^2-\rho_1\rho_2+\rho_2^2-4\pi\rho_1-4\pi\rho_2-2\pi\tau_1-2\pi\tau_2=0.$$
\ethm

We get non-existence of solution for the parameter $\rho$ belonging to two or more regions of the positive quadrant. Such regions are colored in blue in Figure $5$, whereas orange and green regions are the ones for which we have existence of solutions. The pictures show that non-existence phenomena may occur in each of the rectangles where the analysis of Theorem \ref{ex} gave no information. Using the notation of the theorem, these are the cases $(M_1,M_2,M_3)\in\{(1,m,0),(m,1,0),(2,2,1)\}$.\\

We remark that Theorem \ref{nonex2} also applies to the case of $\a_{im}\ge0$. This shows that the existence result in \cite{bjmr} cannot be extended if the hypothesis of positive genus of $\Si$ is removed.\\

\begin{figure}[h!]
\center
\includegraphics[width=0.3\linewidth]{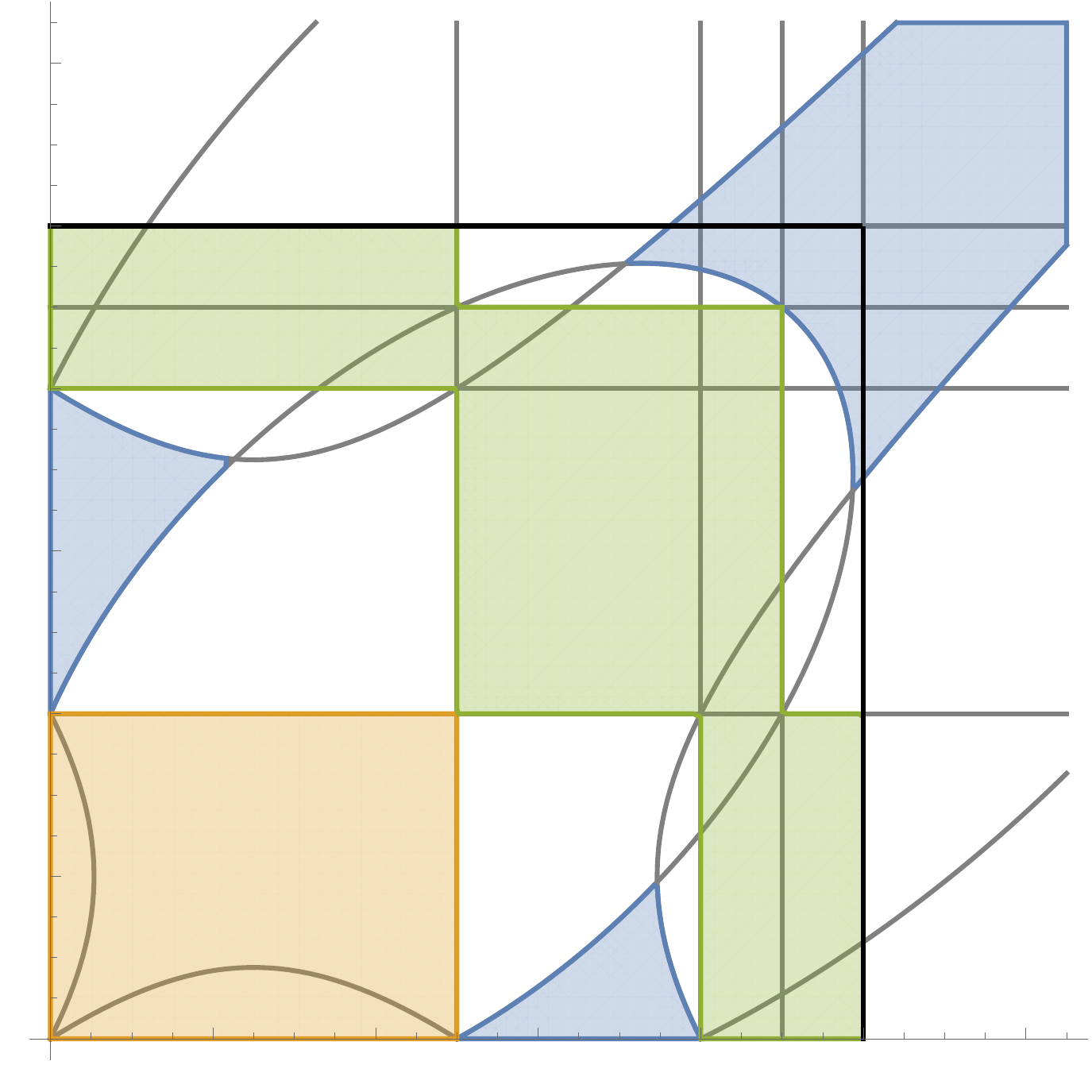}\q\q\q\q\q\q
\includegraphics[width=0.3\linewidth]{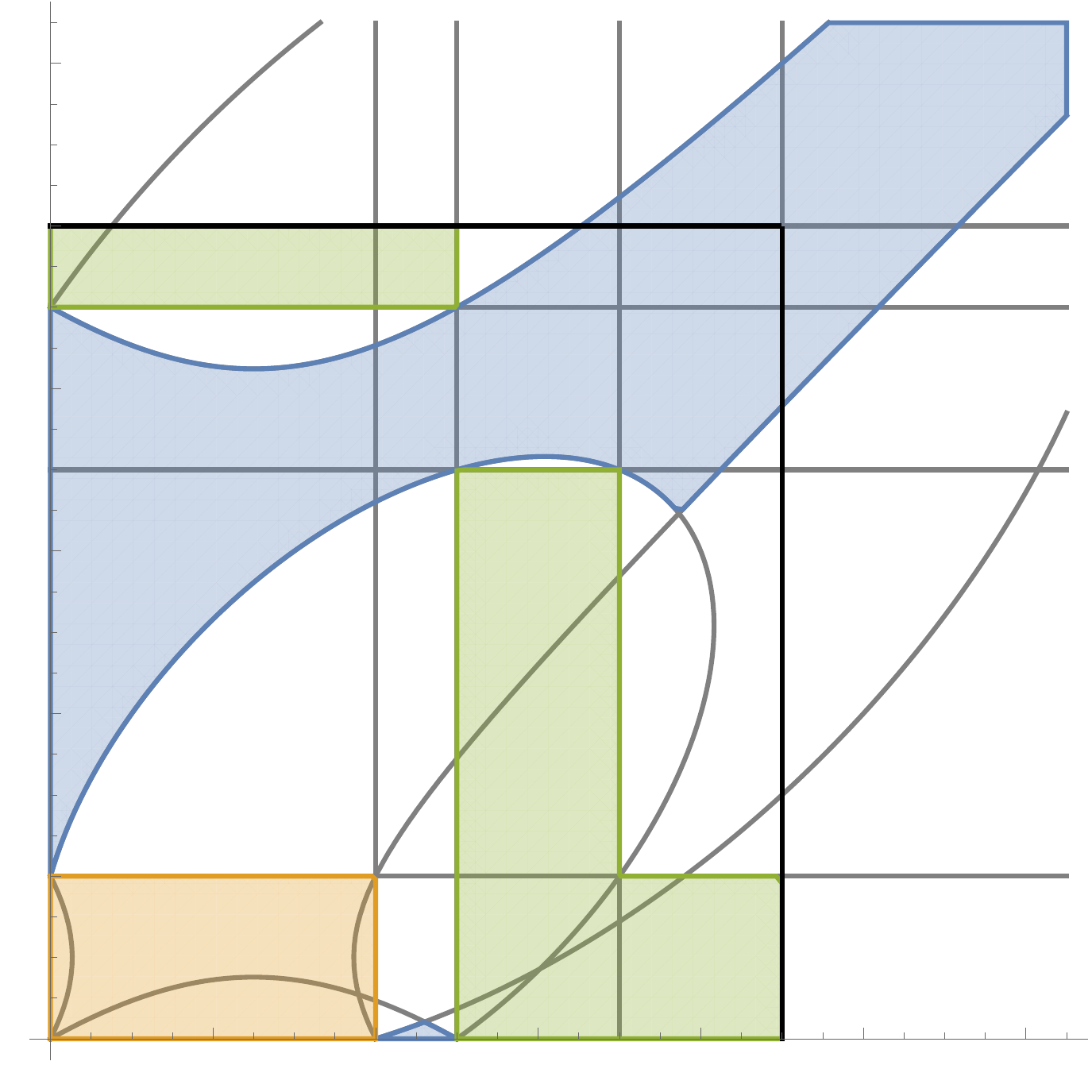}\\
\caption{Values of $\rho$ which yield existence and non-existence results for $\l(\S^2,g_0\r)$, in two different configurations of $\a_{11},\a_{12},\a_{21},\a_{22}$.}
\end{figure}

\bpf[Proof of Theorem \ref{nonex2}]
Let $u=(u_1,u_2)$ be a solution of
$$\l\{\bl-\D_{g_0}u_1=2\rho_1\l(\fr{e^{u_1}}{\int_{\S^2}e^{u_1}dV_{g_0}}-\fr{1}{4\pi}\r)-\rho_2\l(\fr{e^{u_2}}{\int_{\S^2}e^{u_2}dV_{g_0}}-\fr{1}{4\pi}\r)-4\pi\a_{11}\l(\d_{p_1}-\fr{1}{4\pi}\r)-4\pi\a_{12}\l(\d_{p_2}-\fr{1}{4\pi}\r)\\
-\D_{g_0}u_2=2\rho_2\l(\fr{e^{u_2}}{\int_{\S^2}e^{u_2}dV_{g_0}}-\fr{1}{4\pi}\r)-\rho_1\l(\fr{e^{u_1}}{\int_{\S^2}e^{u_1}dV_{g_0}}-\fr{1}{4\pi}\r)-4\pi\a_{22}\l(\d_{p_1}-\fr{1}{4\pi}\r)-4\pi\a_{21}\l(\d_{p_2}-\fr{1}{4\pi}\r)\earr\r.$$
and let $\Pi:\S^2\sm\{p_2\}\to\R^2$ be the stereographic projection. Consider now, for $x\in\R^2$,
$$\l\{\bl U_1(x):=u_1\l(\Pi^{-1}(x)\r)+\log(4\rho_1)-\log\int_{\S^2}e^{u_1}dV_{g_0}-2\a_{11}\log|x|+\l(\fr{\rho_1}{2\pi}-\fr{\rho_2}{4\pi}-\a_{11}-\a_{12}\r)\log\l(1+|x|^2\r)\\
U_2(x):=u_2\l(\Pi^{-1}(x)\r)+\log(4\rho_2)-\log\int_{\S^2}e^{u_2}dV_{g_0}-2\a_{21}\log|x|+\l(\fr{\rho_2}{2\pi}-\fr{\rho_1}{4\pi}-\a_{21}-\a_{22}\r)\log\l(1+|x|^2\r)\earr\r..$$
$U=(U_1,U_2)$ solves
$$\l\{\bl-\D U_1=2H_1e^{U_1}-H_2e^{U_2}\\-\D U_2=2H_2e^{U_2}-H_1e^{U_1}\\\int_{\R^2}H_1(x)e^{U_1(x)}\mrm dx=\rho_1\\\int_{\R^2}H_2(x)e^{U_2(x)}\mrm dx=\rho_2\earr\r.\q\q\q\tx{with}\q\q\q\l\{\bl H_1(x):=\fr{|x|^{2\a_{11}}}{\l(1+|x|^2\r)^{2+\a_{11}+\a_{12}-\fr{\rho_1}{2\pi}+\fr{\rho_2}{4\pi}}}\\H_2(x):=\fr{|x|^{2\a_{21}}}{\l(1+|x|^2\r)^{2+\a_{21}+\a_{22}-\fr{\rho_2}{2\pi}+\fr{\rho_1}{4\pi}}}\earr\r..$$
From Theorem \ref{poho2}, a necessary condition for existence of solutions is
\bequ\label{rhotau}
\rho_1^2-\rho_1\rho_2+\rho_2^2-4\pi\rho_1-4\pi\rho_2-2\pi\tau_1-2\pi\tau_2=0.
\eequ
with $\tau_1,\tau_2$ as in the lemma. Moreover, by the definition of $H_1,H_2$, we have
$$x\cd\n H_i(x)=2\a_{i1}H_i(x)-2\l(2+\a_{i1}+\a_{i2}-\fr{\rho_i}{2\pi}+\fr{\rho_{3-i}}{4\pi}\r)\fr{|x|^2}{1+|x|^2}H_i(x)$$
for both $i$'s, hence we get
$$\tau_i=2\a_{i1}\rho_i-2\l(2+\a_{i1}+\a_{i2}-\fr{\rho_i}{2\pi}+\fr{\rho_{3-i}}{4\pi}\r)\tau'_i\q\q\q\tx{with}\q\q\q\tau'_i:=\int_{\R^2}\fr{|x|^2}{1+|x|^2}H_i(x)\mrm dx.$$
Therefore the necessary condition \eqref{rhotau} becomes
\beqa
\label{taurho}\rho_1^2+\rho_2^2-\rho_1\rho_2-4\pi(1+\a_{11})\rho_1-4\pi(1+\a_{21})\rho_2&+&4\pi\l(2+\a_{11}+\a_{12}-\fr{\rho_1}{2\pi}+\fr{\rho_2}{4\pi}\r)\tau'_1\\
\nonumber&+&4\pi\l(2+\a_{21}+\a_{22}-\fr{\rho_2}{2\pi}+\fr{\rho_1}{4\pi}\r)\tau'_2=0.
\eeqa
Using the straightforward inequalities $0<\tau'_i<\rho_i$ and discussing the cases $2+\a_{i1}+\a_{i2}\lesseqqgtr\fr{\rho_i}{2\pi}-\fr{\rho_{3-i}}{4\pi}$, one can easily see by algebraic computation that \eqref{cond1} and their opposite inequalities are in contradiction with the aforementioned necessary condition.\\
Notice that if $\l\{\bl2+\a_{11}+\a_{12}=\fr{\rho_1}{2\pi}-\fr{\rho_2}{4\pi}\\2+\a_{21}+\a_{22}=\fr{\rho_2}{2\pi}-\fr{\rho_1}{4\pi}\earr\r.$, then \eqref{taurho} just becomes $\rho_1^2+\rho_2^2-\rho_1\rho_2-4\pi(1+\a_{11})\rho_1-4\pi(1+\a_{21})\rho_2=0$. Anyway, one can easily see that these two conditions are equivalent to having all equalities in \eqref{cond1}; this is the reason why we need to assume at least one inequality to be strict.
\epf

\subsection{Proof of Theorems \ref{nonex3} and \ref{nonex4}}

We start by proving Theorem \ref{nonex3}. We will argue by contradiction, following \cite{car} (Theorem $4.1$). Basically, we will assume that a solution exists for some $\a_{11}^n,\a_{12}^n\us{n\to+\infty}\lto-1$. We will consider such a sequence of solutions $u^n$, we will perform a blow-up analysis, following Theorem \ref{conccomp} and we will reach a contradiction.\\

\bpf[Proof of Theorem \ref{nonex3}]
Assume the thesis is false. Then, for some given $\ul\a_{1\wh1},\ul\a_{2\wh1},\rho\nin\G_{\ul\a_{1\wh1},\ul\a_{2\wh1}}$, there exist a sequence $(\a_{11}^n,\a_{21}^n)\us{n\to+\infty}\lto(-1,-1)$ and a sequence $u^n=(u_1^n,u_2^n)$ of solutions of
$$\l\{\bl-\D u_1^n=2\rho_1\l(\fr{\wt h_1^ne^{u_1^n}}{\int_\Si\wt h_1^ne^{u_1^n}\mrm dV_g}-1\r)-\rho_2\l(\fr{\wt h_2^ne^{u_2^n}}{\int_\Si\wt h_2^ne^{u_2^n}\mrm dV_g}-1\r)\\-\D u_2^n=2\rho_2\l(\fr{\wt h_2^ne^{u_2^n}}{\int_\Si\wt h_2^ne^{u_2^n}\mrm dV_g}-1\r)-\rho_1\l(\fr{\wt h_1^ne^{u_1^n}}{\int_\Si\wt h_1^ne^{u_1^n}\mrm dV_g}-1\r)\earr\r.,$$
with $\wt h_1^n,\wt h_2^n$ such that $\wt h_i^n\sim d(\cd,p_1)^{2\a_{i1}^n}$. It is not restrictive to assume
$$\int_\Si\wt h_1^ne^{u_1^n}\mrm dV_g=\int_\Si\wt h_2^ne^{u_2^n}\mrm dV_g=1.$$
We would like to apply Theorem \ref{conccomp} to the sequence $u^n$. Anyway, since the coefficients $\a_{i1}^n$ are not bounded away from $-1$, we cannot use such a theorem on the whole $\Si$, but we have to remove a neighborhood of $p_1$. A first piece of information about blow-up is given by the following Lemma, inspired by \cite{car}, Lemma $4.3$.
\epf

\blem\label{blowup}
Let $\d>0$ small be given and $u^n$ be as in the proof of Theorem \ref{nonex3}. Then, $u_1^n,u_2^n$ cannot be both uniformly bounded from below on $\pa B_\d(p_1)$.
\elem

\bpf
Assume by contradiction that $\inf_{\pa B_\d(p_1)}u_i^n\ge-C$ for both $i$'s and define $v^n:=\fr{2u_1^n+u_2^n}3$. Then
$$\l\{\bll-\D v^n=\rho_1\l(\wt h_1^ne^{u_1^n}-1\r)\ge-\rho_1&\tx{in }B_\d(p_1)\\v^n\ge-C&\tx{on }\pa B_\d(p_1)\earr\r..$$
By the maximum principle, $v^n\ge-C$ on $B_\d(p_1)$, therefore by the convexity of the exponential function we get the following contradiction:
\beqy
+\infty&\us{n\to+\infty}\lot&\int_{B_\d(p_1)}d(\cd,p_1)^{2\max\{\a_{11}^n,\a_{21}^n\}}\mrm dV_g\\
&\le&C\int_{B_\d(p_1)}d(\cd,p_1)^{2\max\{\a_{11}^n,\a_{21}^n\}}e^{v^n}\mrm dV_g\\
&\le&C\l(\fr{2}3\int_{B_\d(p_1)}d(\cd,p_1)^{2\max\{\a_{11}^n,\a_{21}^n\}}e^{u_1^n}\mrm dV_g+\fr{1}3\int_{B_\d(p_1)}d(\cd,p_1)^{2\max\{\a_{11}^n,\a_{21}^n\}}e^{u_2^n}\mrm dV_g\r)\\
&\le&C\l(\int_{B_\d(p_1)}\wt h_1^ne^{u_1^n}\mrm dV_g+\int_{B_\d(p_1)}\wt h_2^ne^{u_2^n}\mrm dV_g\r)\\
&\le&C.
\eeqy
This concludes the proof.
\epf

\bpf[Proof of Theorem \ref{nonex3}, continued]
Let us apply Theorem \ref{conccomp} to $u^n$ on $\O:=\Si\sm B_\fr{\d}2(p_1)$ for some given small $\d>0$. By Lemma \ref{blowup}, boundedness from below cannot occur for both components, therefore we either have blow-up or (up to switching the indexes) $u_1^n\us{n\to+\infty}\lto-\infty$ uniformly on $\Si\sm B_\d(p_1)$. In other words,
$$\left.\rho_1\fr{\wt h_1^ne^{u_1^n}}{\int_\Si\wt h_1^ne^{u_1^n}\mrm dV_g}\right\lfloor_{\Si\sm B_\d(p_1)}\us{n\to+\infty}\wk\sum_{x\in\mcal S}\s_1(x)\d_x\q\q\left.\rho_2\fr{\wt h_2^ne^{u_2^n}}{\int_\Si\wt h_2^ne^{u_2^n}\mrm dV_g}\right\lfloor_{\Si\sm B_\d(p_1)}\us{n\to+\infty}\wk r_2+\sum_{x\in\mcal S}\s_2(x)\d_x,$$
where we set $\mcal S=\es$ if blow up does not occur. Anyway, being $\d$ arbitrary, a diagonal argument gives
$$\rho_1\fr{\wt h_1^ne^{u_1^n}}{\int_\Si\wt h_1^ne^{u_1^n}\mrm dV_g}\us{n\to+\infty}\wk\sum_{x\in\mcal S}\s_1(x)\d_x+\s_1(p_1)\d_{p_1}\q\q\rho_2\fr{\wt h_2^ne^{u_2^n}}{\int_\Si\wt h_2^ne^{u_2^n}\mrm dV_g}\us{n\to+\infty}\wk r_2+\sum_{x\in\mcal S}\s_2(x)\d_x+\s_2(p_1)\d_{p_1},$$
with $\s_1(p_1)=\rho_1-\sum_{x\in\mcal S}\s_1(x)$ and $\s_2(p_1)=\rho_2-\sum_{x\in\mcal S}\s_2(x)-\int_\Si r_2\mrm dV_g$.\\
By a variation of the Poho\v zaev identity (see \cite{lwzhang}, Proposition $3.1$ and \cite{batman}, Lemma $2.4$), we get
$$\s_1(p_1)^2-\s_1(p_1)\s_2(p_1)+\s_2(p_1)^2=0,$$
that is $\s_1(p)=\s_2(p)=0$. In particular, we get $\rho_1=\sum_{x\in\mcal S}\s_1(x)$, which means either $\rho_1=0$ or $\rho\in\G_{\ul\a_{1\wh1},\ul\a_{2\wh1}}$ This contradicts the assumptions and proves the theorem.
\epf

\

We conclude by proving Theorem \ref{nonex4}. The argument is somehow similar: we assume, by contradiction, to have a solution satisfying all the hypotheses for $\a_{23}^n\us{n\to+\infty}\lto-1$. Then, we perform a blow-up analysis and we rule out the last case using Theorem \ref{nonex2}.\\

\bpf[Proof of Theorem \ref{nonex4}]
Assume by contradiction there exists a sequence $\a_{23}^n\us{n\to+\infty}\lto-1$ such that \eqref{toda} has a solution $u^n$ satisfying all the hypotheses of Theorem \ref{nonex4} and, w.l.o.g., $\int_\Si\wt h_1^ne^{u_1^n}\mrm dV_g=\int_\Si\wt h_2^ne^{u_2^n}\mrm dV_g=1$. In particular, since $4\pi(1+\a_{2m})<\rho_2^n<4\pi(2+\a_{2m}+\a_{23}^n)$ for $m=1,2$, then $\rho_2^n\us{n\to+\infty}\lto4\pi(1+\a_{21})=4\pi(1+\a_{22})$.\\
As in the proof of Theorem \ref{nonex3}, we must have $\inf_{B_\d(p_3)}u_2^n\us{n\to+\infty}\lto-\infty$ for small $\d$, because otherwise we would have
$$1\ge\int_{B_\d(p_3)}\wt h_2^ne^{u_2^n}\mrm dV_g\us{n\to+\infty}\ge C\int_{B_\d(p_3)}d(\cd,p_3)^{2\a_{23}^n}\mrm dV_g\us{n\to+\infty}\lto+\infty.$$
Unlike before, we cannot apply the maximum principle to get $\inf_{\pa B_\d(p_3)}u_2^n\us{n\to+\infty}\lto-\infty$. Anyway, this could be ruled out by the following argument: if $u_2^n$ were uniformly bounded from below on $\pa B_\d(p_3)$ and $u_2^n(x^n)\us{n\to+\infty}\lto-\infty$ for some $x^n\us{n\to+\infty}\lto x\in B_\d(p_3)$, then applying Theorem \ref{conccomp} we would get blow-up at $x$ of the first component alone, which would give $\rho_1\ge\s_1(x)\ge4\pi$, in contradiction to the assumptions of the theorem. Therefore, $u_2^n$ must go to $-\infty$ uniformly on $\pa B_\d(p_3)$, which means, by Theorem $\ref{conccomp}$, blow-up with $r_2\eq0$.\\
The assumption $\rho_2<4\pi$ implies that such blow-up must occur at a subset of $\{p_1,p_2,p_3\}$. Blow-up in $p_3$ is also excluded because, by standard blow-up analysis (one can argue for instance as in \cite{os}, Lemma $9$) it would imply $\rho_1\ge\s_1(p_3)\ge4\pi$; therefore, we must have $\wt h_2e^{u_2^n}\us{n\to+\infty}\wk4\pi(1+\a_{2m})\d_{p_m}$ for some $m=1,2$; since the role of $p_1$ and $p_2$ is interchangeable, we can assume $m=1$.\\
Let us now consider $u_1^n$: it cannot blow up at $p_1$, because $\rho_1<4\pi(2+\a_{11}+\a_{12})$, and it cannot even blow up at any other points: in fact, since $u_2^n$ only blows up at $p_1$, then by Theorem \ref{conccomp} we would get $\rho_1=\sum_{x\in\mcal S}\s_1(x)$, again contradicting the assumptions. Therefore, $u_1^n$ must converge to some $u_1$, which solves (up to subtracting a suitable combination of Green's functions)
$$-\D u_1=2\rho_1\l(\fr{e^{u_1}}{\int_{\S^2}e^{u_1}\mrm dV_{g_0}}-\fr{1}{4\pi}\r)-4\pi(2+\a_{11}+\a_{21})\l(\d_{p_1}-\fr{1}{4\pi}\r)-4\pi\a_{12}\l(\d_{p_2}-\fr{1}{4\pi}\r).$$
By applying Theorem \ref{nonex2} with $\rho_2=0$, or equivalently Proposition $5.8$ in \cite{barmal}, we see that the last equation is not solvable since $4\pi(1+\a_{12})<\rho_1<4\pi(2+\a_{11}+\a_{21})$. This gives a contradictions and proves the theorem.
\epf

\

Notice that, by repeating the same argument, we can find similar non-existence results in the cases $(M_1,M_2,M_3)=(2,2,1)$ and $(1,m,0)$, with one coefficient being very close to $-1$. In the case $(1,m,0)$, we can also drop the assumptions on $\Si$ to be a standard sphere with antipodal singularities.

\bibliography{finale}
\bibliographystyle{abbrv}

\end{document}